\documentclass{article}
\providecommand{\keywords}[1]{\textbf{Keywords:} #1}
\usepackage{lineno,hyperref}
\usepackage{bm,upgreek}
\usepackage[a4paper, left=2cm, right=2cm, top=2cm]{geometry}
\usepackage{color}
\usepackage{changes}
\usepackage{graphicx}
\usepackage[english]{babel}
\usepackage{float}
\usepackage{amsmath}
\usepackage{amssymb}
\usepackage{tabularx}
\usepackage{booktabs}
\usepackage{multirow}
\usepackage{todonotes}
\usepackage{makecell}
\usepackage{authblk}
\usepackage{fancyhdr}

\usepackage[backend=biber,citestyle=numeric-comp,bibstyle=ieee,dashed=false,url=false]{biblatex}
\AtEveryBibitem{\ifentrytype{article}{\clearfield{note}}{}}
\AtEveryBibitem{\ifentrytype{book}{\clearfield{note}}{}}
\AtEveryBibitem{\ifentrytype{article}{\clearfield{issn}}{}}
\DeclareCiteCommand{\cite}
    {\usebibmacro{cite:init}%
     \usebibmacro{prenote}}
    {\usebibmacro{citeindex}%
         \ifnumequal{\value{citecount}}{1}{\printtext[bibhyperref]{[}}{}%
        \usebibmacro{cite:comp}%
    }
    {}
    {\usebibmacro{cite:dump}%
     \ifnumequal{\value{citecount}}{\value{citetotal}}{\printtext[bibhyperref]{]}}{}%
     \usebibmacro{postnote}%
    }
\begin{document}

\title{Nonlinear synthesis of compliant mechanisms with selective compliance}

\author{Stephanie Seltmann}

\author{Alexander Hasse*}

\affil{Chemnitz University of Technology, Professorship Machine Elements and Product Development, Chemnitz, 09107, Germany}

\maketitle

\begin{abstract}
The synthesis of compliant mechanisms (CMs) is frequently achieved through topology optimization. Many synthesis approaches simplify implementation by assuming small distortions, but this limits their practical application since CMs typically undergo large deformations that include geometric and material nonlinearities. CMs designed to generate a desired deformation path at the output points under specific loads are known as path-generating CMs. However, these CMs face significant challenges in topology optimization, resulting in the development of only a few optimization methods. Existing approaches often include only certain load cases in the optimization process. Consequently, if a CM designed this way encounters different load cases in practice, its path-generating behavior cannot be guaranteed.

The authors have previously contributed to the development of an approach suitable for synthesizing load case insensitive CMs. This paper extends that approach to account for nonlinearities, enabling the synthesis of path-generating CMs. The effectiveness of this extended approach is demonstrated through appropriate design examples. Additionally, the paper presents, for the first time, a shape-adaptive path-generating CM.
    
\vspace{0.5cm}
\noindent

\end{abstract}
\fancypagestyle{firstpage}{
  \fancyhf{}
  \fancyfoot[L]{* Corresponding author.
  E-mail address: alexander.hasse@mb.tu-chemnitz.de}
  \renewcommand{\headrulewidth}{0pt}
  \renewcommand{\footrulewidth}{0.4pt}
}

\thispagestyle{firstpage} 

\section{Introduction}
\label{sec:Introduction}

\begin{figure}
    \centering
    \includegraphics[width=0.8\textwidth]{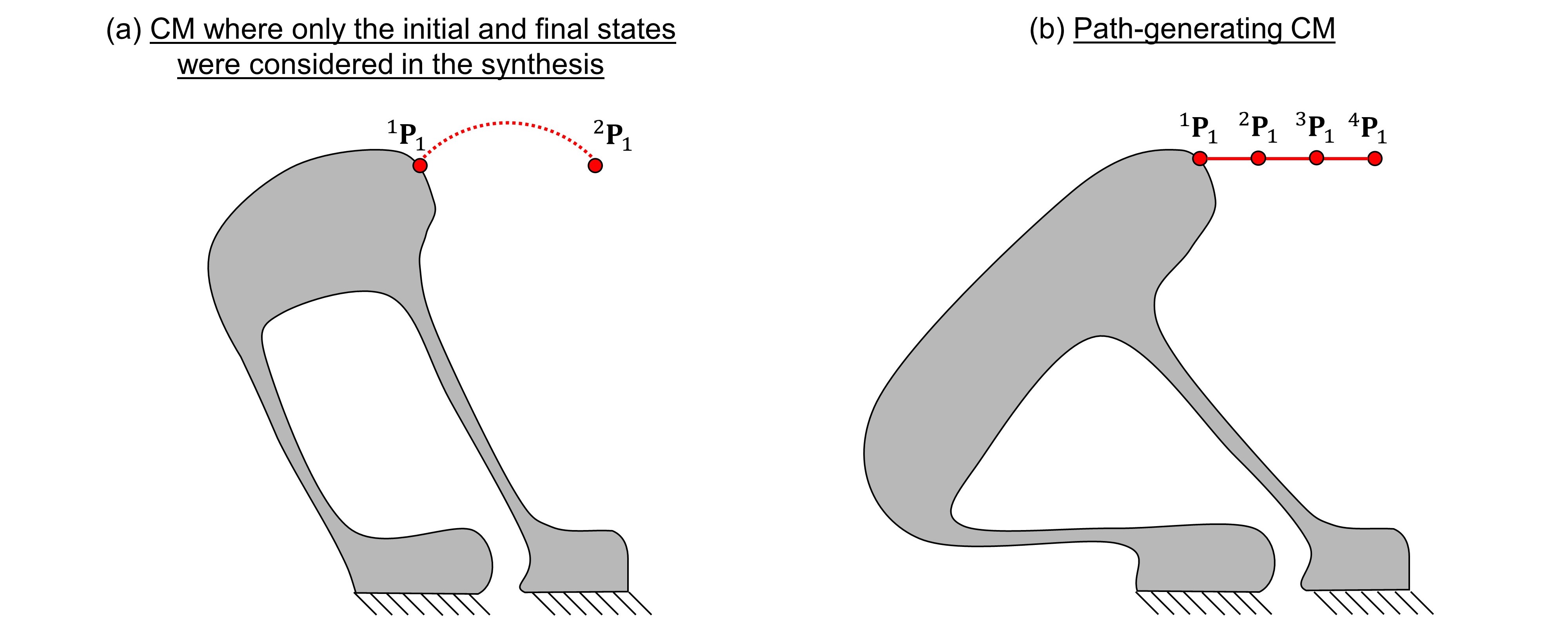}
    \caption{Types of nonlinear CM \label{fig0}}
\end{figure}

The topology optimization of compliant mechanisms (CMs) typically relies on linear assumptions, including small deformations and linear-elastic material behaviour. However, these assumptions do not accurately reflect real-world conditions, wherein CMs often experience large deformations. Consequently, it is imperative to account for nonlinearities, encompassing both geometric nonlinearity (large deformations) and material nonlinearity, to achieve more precise and effective designs \cite{Bruns1998,DeLeon2020}.
CMs are typically designed to ensure a specific deformation of the output points when the input points are subjected to certain actuator loads or predefined displacements \cite{Deepak2009, Bendsoe2003d}. In addition to the requirement for the outputs to achieve the desired deformations, most optimization approaches also apply transverse loads to these points to attain a particular stiffness in the CM. However, if different load cases are encountered in practice, it may no longer be feasible to guarantee the desired deformation.
In contrast, the kinematic design of conventional mechanisms does not require the specification of loads, as the synthesis is purely driven by kinematics. The authors have contributed to the development of the so-called pseudo-kinematic approaches \cite{Hasse2009, Hasse2016, Kirmse2021}, which adopt a similar strategy for the design of CMs. In this approaches, the synthesis of a CM is typically based on a \textit{desired kinematics}, meaning the specific kinematics that the CM should exhibit during actuation. This desired kinematics thereby serves as main design parameter – no reference to loads is made. However, the pseudo-kinematic approaches are based on linear assumptions in their current form. This paper focuses on extending the approaches to incorporate the aforementioned nonlinearities.
Before delving into the specifics of the pseudo-kinematic approaches, this paper provides a brief overview of existing nonlinear topology optimization methods, highlighting the associated challenges and research gaps. The authors distinguish between different types of nonlinear synthesis methods for CM. The first type considers only the initial and final deformation states (see e.g. \cite{Bruns1998, Pedersen2001, DeLeon2020, Joo2004, Maute2003, Li2014a, Lazarov2011, Capasso2020}). In these approaches, the deformation path between the initial state $^1\mathbf{P}$ and the final state $^2\mathbf{P}$ is not determined during optimization, resulting in arbitrary intermediate deformation states. Figure~\ref{fig0}a illustrates this for a CM with a single output point $\mathbf{P}_1$, where the degrees of freedom (DoFs) of this point are the output DoFs. These methods are typically limited to CMs with a few output DoFs. In contrast, the approach presented in \cite{Kumar2021} is also suitable for the synthesis of shape-adaptive structures, which are characterized by having many output DoFs and undergoing surface deformation as part of the desired deformation.

\unskip

The second type of approach considers the deformation path of the compliant mechanism between the initial and final states during optimization, using stationary points. These resulting CMs are referred to as path-generating CMs (see Figure~\ref{fig0}b). However, according to \cite{Zhu2020}, this type of CM presents significant challenges in topology optimization: "Topology optimization of path generating CMs is a highly complex topology optimization problem. It requires a considerable amount of experiments in formulating objective functions as well as nonlinear finite element analysis". Consequently, only a few optimization methods have been developed for this purpose. \textsc{Pedersen et al.} \cite{Pedersen2001} propose an objective function that minimizes the deviation of the output DoFs' displacement between the target and actual deformation paths. This deviation is measured at specific stationary points (for defined input displacements) and summarized. Additionally, counter loads are introduced at these stationary points as transverse loads to ensure meaningful results. In \cite{Saxena2001a}, this deviation function is also minimized. Some examples in this study do not include transverse loads, while one example incorporates transverse loads caused by springs. \textsc{Mankame et al.} \cite{Mankame2007} introduce an objective function based on Fourier descriptors, where the deviation between the target and actual deformation paths described by Fourier descriptors is minimized. The optimization process includes certain transverse load cases at the output.
The approaches presented thus far for path-generating CMs predominantly produce CMs with lumped compliance. These mechanisms suffer from the disadvantage of stress concentrations in the compliant areas, which limits their damage-free deformability. Alternatively, compliant mechanisms with distributed compliance, as described by \cite{Kota1995a}, do not have this drawback. By distributing the compliant areas throughout the entire CM, stress concentrations are avoided. However, synthesizing such CMs requires further development of existing optimization formulations \cite{Seltmann2023}. An example of this is found in \cite{daSilva2020}, which extends the objective function used in \cite{Pedersen2001} by incorporating stress constraints and production inaccuracies, thereby creating CMs with distributed compliance. Nonetheless, these stress constraints render the optimization problem highly non-convex, making it challenging to satisfy them precisely \cite{daSilva2021}.

In summary, the review of synthesis methods for compliant mechanisms reveals several shortcomings. Firstly, the designed CMs typically have only one output point, highlighting a lack of approaches for synthesizing shape-adaptive, path-generating CMs. Secondly, only few methods are suitable for designing CMs with distributed compliance. Thirdly, current approaches primarily focus on minimizing the deviation between the desired and actual deformation paths, referred to as accuracy in the literature \cite{Campanile2022}. Additionally, the deviations in the deformation path under different transverse loads, known as precision \cite{Campanile2022}, are insufficiently considered.
An approach is therefore required to address all three identified shortcomings. Under linear assumptions, pseudo-kinematic approaches have been developed to comprehensively consider both accuracy and precision during synthesis by optimizing the eigenbehavior of the structure \cite{Hasse2009, Hasse2016, Kirmse2021}. This involves minimizing the ratio between the first and second eigenvalues to increase selectivity, thereby aligning the first eigenvector with the desired kinematics. Selectivity measures the precision of the mechanism. Further explanations on this can be found in section 5. Moreover, pseudo-kinematic approaches are suitable for synthesizing CMs with multiple output DoFs. Additionally, the method in \cite{Kirmse2021} was extended in \cite{Seltmann2023} to produce CMs with selective compliance, a subtype of distributed compliance that exhibits high precision. This method is now being extended for application to nonlinear path-generating CMs. To achieve this, it will be applied to a set of stationary points corresponding to the desired kinematics in the output DoFs. The desired properties can be realized by optimizing the eigenbehavior at these stationary points. This new optimization of eigenbehavior for the synthesis of path-generating CMs fundamentally differs from the deviation-based optimization formulations used in other methods.

\section{Path-generating compliant mechanisms with selective compliance}
\label{Pathgenerating_CM}

Since they realize the deformation through elastic distortion, CMs must be modelled as continua. This means that they have an infinite number of structural DoFs. However, for analysis and synthesis using FEM, these DoFs are discretized to a finite number $p$. As already mentioned, other authors divide DoFs, which either perform certain predefined displacements or are subject to external forces, into input and output DoFs. In the synthesis of CM with selective compliance, this subdivision is omitted and the DoFs are summarized as $q$ active DoFs and labelled with the subscript $a$. The remaining structural DoFs are labelled as $p-q$ passive DoFs (subscript $c$).

As aforementioned, the pseudo-kinematic approaches for the synthesis of CM with selective compliance have so far been based on linear assumptions \cite{Hasse2009,Hasse2016,Kirmse2021,Seltmann2022,Seltmann2023}. These assumptions are used to establish and solve the balance of internal and external forces on the structure to be designed with respect to the unloaded configuration. This means that the stiffness matrix $\mathbf{K}_\mathrm{L}$ of the CM can be calculated independently of the deformations, which greatly simplifies the Finite Element Method (FEM). The characteristic feature of the CMs with selective compliance designed to date is that, regardless of the loads acting on the active DoFs, a displacement occurs that corresponds with a small deviation to a scaling $\alpha_\mathrm{L}$ (with the unit of a length) of a desired deformation mode $\bar{\bm{\upvarphi}}$:

\begin{linenomath}
    \begin{equation}
        \label{desired_motion_linear}
        \mathbf{u}_a\approx\upalpha_\mathrm{L}\bar{\bm{\upvarphi}}
    \end{equation}
\end{linenomath}

The desired deformation mode is specified as input parameter for the optimization and defines a specific displacement pattern that the active DoFs should perform. Deformations corresponding to the desired kinematics are also referred as \textit{desired deformations}. Due to the previous limitations, the deformation behavior of the CM can only be defined for small deformations. They are referred to below as linear CM.

In order to design path-generating CMs, the restriction to small distortions has to be eliminated. In this case, the stiffness matrix $\mathbf{K}$ of the CM is no longer constant, but changes depending on the deformations of the CM. The previous form \eqref{desired_motion_linear} can no longer be used to describe the desired deformations. The desired deformations are now represented by a non-linear vector function, the \textit{desired deformation function} $\bar{\bm{\upxi}}(\upalpha)$. In the case of a path-generating CM with selective compliance, a deformation in the active DoFs is ideally element of the desired deformation function independently of the applied load:

\begin{equation}
    \label{desired_motion_nonlinear}
    \begin{bmatrix}
    u_{a1}\\u_{a2}\\\vdots\\u_{aq}
    \end{bmatrix}
    \approx 
    \begin{bmatrix}
    \bar{\upxi}_1(\upalpha)\\\bar{\upxi}_2(\upalpha)\\\vdots\\\bar{\upxi}_q(\upalpha)
    \end{bmatrix}
    \iff
    \mathbf{u}_{a}\approx \bar{\bm{\upxi}}(\upalpha)\mid\upalpha\in[\upalpha_{\mathrm{l}},\upalpha_{\mathrm{u}}]
\end{equation}

The desired deformations are described in a certain domain of definition for $\upalpha$, which is limited by the boundaries $\upalpha_{\mathrm{l}}$ and $\upalpha_{\mathrm{u}}$. The component functions $\bar{\upxi}(\upalpha)$ have to be continuous and differentiable functions whose derivative with respect to $\upalpha$ does not become $\mathbf{0}$ in the domain of definition. In addition, the desired deformation function has to include the undeformed state. The component functions in the parameter representation therefore do not contain any constant components. This requirement is based on the assumption that the CMs are not prestressed. The vector function provides a vector with the unit of a length as the function value. The newly introduced desired deformation function can be used to describe both linear and non-linear deformations.

The desired deformation function represents the desired kinematics in continuous form. However, this continuous form cannot be used algorithmically when using the FEM for structural modelling during synthesis. It must therefore be converted into a discrete form. Therefore, the desired deformation function is discretized to a number of $n$ stationary points $t$ and linearized in these stationary points.

To determine the stationary points, values for $\upalpha$ are specified and then the associated function values are determined, which are subsequently referred to as $^t\mathbf{u}_{a}$. In addition, the derivatives of the desired deformation function $\bar{\bm{\upxi}}(\upalpha)$ are determined at the stationary points, which still have the unit of a length. Normalization makes the derivatives dimensionless and results in the desired tangent deformation modes $^t\bar{\bm{\upvarphi}}$. For small deviations from the stationary points, the desired deformation can be linearized by scaling the desired tangent deformation mode $^t\bar{\bm{\upvarphi}}$ belonging to the stationary point similarly to \eqref{desired_motion_linear}. The linearized equation is

\begin{equation}
    \label{desired_motion_linearized}
    ^t\mathbf{u}_a+\Delta\mathbf{u}_a\approx~{^t\mathbf{u}_a}+\alpha_\mathrm{L}~{^t\bar{\bm{\upvarphi}}},\;t=1\ldots n
\end{equation}

The tangent stiffness matrices in the stationary points $\mathbf{K}_\mathrm{NL}({^t\mathbf{u}}_a)$ have to be known as a precondition for the optimization of a path-generating CM with selective compliance. These are a linearization of the stiffness matrix $\mathbf{K}_\mathrm{NL}(\mathbf{u})$ at the deformation $^t\mathbf{u}_a$ in the current stationary point. They can therefore be regarded as constant for small deviations from this deformation and calculations with the linearized stiffness matrices can be carried out equivalently to linear assumptions.
\section{Equilibrium analysis}
\label{Equilibrium analysis}

In the case of linear assumptions (linear FEM), the stiffness matrix $\mathbf{K}_\mathrm{L}$ and the internal nodal forces $\mathbf{f}_\mathrm{L}$ of a set of $m$ elements can be calculated as follows:

\begin{linenomath}
    \begin{equation}
        \label{linear_K}
        \mathbf{K}_\mathrm{L}=\sum_{e=1}^m\int_{V^{e}}\mathbf{B}_\mathrm{L0}^{e\mathrm{T}}~\mathbf{C}_\mathrm{L}^{e}~\mathbf{B}_\mathrm{L0}^{e}~dV^{e}
    \end{equation}
\end{linenomath}

\begin{linenomath}
    \begin{equation}
        \label{linear_F}
        \mathbf{f}_\mathrm{L}=\mathbf{K}_\mathrm{L}\mathbf{u}
    \end{equation}
\end{linenomath}

Here, $\mathbf{C}_\mathrm{L}^{e}$ is the constant stress-strain material property matrix for the element $e$. In the following, a plane stress state and isotropic material are assumed. $\mathbf{B}_\mathrm{L0}^{e}$ is the constant distortion-displacement matrix.

For the non-linear FEM, the following equations are based on the total Lagrangian formulation. The tangent stiffness matrix and the internal nodal forces at the stationary point $t$ can be calculated as follows (for a detailed derivation see \cite{Bathe2014}):

\begin{linenomath}
    \begin{equation}
        \label{Nonlinear_K}
        ^t\mathbf{K}_\mathrm{NL}=\sum_{e=1}^m\int_{V^{e}}{^t\mathbf{B}^{e\mathrm{T}}_\mathrm{L}}~ {^t\mathbf{C}^{e}}~{^t\mathbf{B}^{e}_\mathrm{L}}~dV^{e}+\sum_{e=1}^m\int_{V^{e}} {^t\mathbf{B}^{e\mathrm{T}}_\mathrm{NL}}~{^t\mathbf{S}^{e}}~{^t\mathbf{B}^{e}_\mathrm{NL}}~dV^{e}
    \end{equation}
\end{linenomath}

\begin{linenomath}
    \begin{equation}
        \label{Nonlinear_F}
        ^t\mathbf{f}_\mathrm{NL}=\sum_{e=1}^m\int_{V^{e}} {^t\mathbf{B}^{e\mathrm{T}}_\mathrm{L}}~ {^t\mathbf{\hat{S}}^{e}}~dV^{e}
    \end{equation}
\end{linenomath}

This formulation includes the effects of large displacements and rotations (geometric non-linearities). The strain-displacement transformation matrices ${^t\mathbf{B}^{e}_\mathrm{L}}$ and ${^t\mathbf{B}^{e}_\mathrm{NL}}$ are dependent on $\mathbf{u}$. The stress-strain material property matrix ${^t\mathbf{C}^{e}}$ can either be assumed to be constant for the case of linear elasticity ($\mathbf{C}_\mathrm{L}^{e}$) or dependent on the displacements ($^t\mathbf{C}_\mathrm{NL}^{e}$) if material non-linearity is taken into account. ${^t\mathbf{S}^{e}}$ and ${^t\mathbf{\hat{S}}^{e}}$ are the second Piola-Kirchhoff stresses, but the same entries are arranged at different positions. 4-node quadrilateral elements are used for the topology optimization. The volume integrals contained in equation \eqref{linear_K} to \eqref{Nonlinear_F} are calculated using full Gauss numerical integration.

The finite element equation can only be calculated approximately for non-linear FEM, as both $^t\mathbf{K}_\mathrm{NL}$ and $^t\mathbf{f}_\mathrm{NL}$ depend on $\mathbf{u}$. The structural equilibrium can be written as follows:

\begin{linenomath}
    \begin{equation}
        \label{Residuum}
        ^t\mathbf{r} - {^t\mathbf{f}_\mathrm{NL}} =\mathbf{0}
    \end{equation}
\end{linenomath}

where $^t\mathbf{r}$ are the external forces. This equilibrium is determined using the Newton-Raphson method. For a detailed derivation, please refer to \cite{Bathe2014}. If the difference between external and internal nodal forces is smaller than a convergence tolerance $\upepsilon_\mathrm{F}$, the Newton-Raphson iteration is finished.

\section{Dealing with numerical difficulties in structural modelling}
\label{Numerical_Problems}

Structural modelling with non-linear FEM causes numerical difficulties in the context of topology optimization. Therefore, the equilibrium analysis must first be stabilized for the subsequent topology optimization.

At the beginning, the stiffness matrices for all stationary points are parameterized. This is done by scaling the stiffnesses of the element stiffness matrices $^t\mathbf{K}^e$ for each stationary point with a design variable $x^{e}$ before they are integrated into the stiffness matrix:

\begin{linenomath}
    \begin{equation}
        \label{Parametrization}
        ^t\mathbf{K}(\mathbf{x})= \sum_{e=1}^m x^{e}~{^t\mathbf{K}^e}, \; t=1...n    
    \end{equation}
\end{linenomath}

The matrix $^t\mathbf{K}^e$ is a stiffness matrix that depends on $^t\mathbf{u}$, which means that there is a different $^t\mathbf{K}$ for each stationary point. The matrices $^t\mathbf{K}$ for each stationary point are all parameterized with the same set of design variables. The matrix $^t\mathbf{K}^e$ does not correspond to $^t\mathbf{K}_\mathrm{NL}$; the necessary adjustments will be explained below. In contrast to the well-known Solid Isotropic Material with Penalization (SIMP) approach \cite{Bendsoe2003c}, $x^{e}$ is used for the parameterization instead of $(x^{e})^{\upeta}$, where $\upeta$ is a penalty factor. Further explanations can be found in the section \ref{Problem_statement}. During optimization, a constraint is usually introduced whereby the design variables can only take values between a small positive number and one. 

If the associated design variable of an element is small, the element has a low stiffness and tends to excessive distortion when the CM is deformed \cite{Bruns2003}. These low stiffness elements hinder the convergence of the Newton-Raphson iterations. The convergence is strongly influenced by the material model used. If the Saint Venant-Kirchhoff model is used, geometric non-linearities are included by using the second Piola-Kirchhoff stress tensor, but the stress-strain material property matrix ${^t\mathbf{C}^{e}}$ is assumed to be constant ($\mathbf{C}_\mathrm{L}^{e}$). This model is only suitable for modelling small distortions; for larger distortions, the finite elements show unrealistic deformation behavior. However, during optimization, small distortions cannot be ensured in all finite elements \cite{Klarbring2013}. The use of a polyconvex hyperelastic material model can improve the convergence, as it has a stiffening effect during compressive deformations in the finite elements and enables a more realistic deformation behavior with larger distortions. This is shown for example in \cite{Bruns1998, DeLeon2020, Lahuerta2013}. For this reason, we also use a hyperelastic material model with the strain energy function of the compressible neo-Hookean material model presented in \cite{Simo1998a}:

\begin{linenomath}
    \begin{equation}
        \label{strain-energy_function}
        ^t\psi_\mathrm{NL}=\uplambda_\mathrm{M}(\frac{{^tJ^2}-1}{4})-(\frac{\uplambda_\mathrm{M}}{2}+\upmu_\mathrm{M})\ln ({^tJ})+\frac{1}{2}\upmu_\mathrm{M}(^tC_{kk}-3)
    \end{equation}
\end{linenomath}

Here $^tC_{ij}$ is the right Cauchy-Green deformation tensor and $^tJ^2$ its determinant. The Lamé parameters $\upmu_\mathrm{M}$ and $\uplambda_\mathrm{M}$ with $E_\mathrm{M}$ as the modulus of elasticity and $\upnu_\mathrm{M}$ as the Poisson's ratio and the Kroneker delta are defined as usual:

\begin{linenomath}
    \begin{equation}
        \label{Lame_parameters}
        \uplambda_\mathrm{M}=\frac{E_\mathrm{M}\upnu_\mathrm{M}}{(1+\upnu_\mathrm{M})(1-2\upnu_\mathrm{M})}; \; \upmu_\mathrm{M}=\frac{E_\mathrm{M}}{2(1+\upnu_\mathrm{M})};
    \end{equation}
\end{linenomath}

\begin{linenomath}
    \begin{equation}
        \label{Kroneker_delta}
        \updelta_{ij}=\left\{
        \begin{matrix}
            0; \; i \neq j \\ 1; \; i=j
        \end{matrix}
        \right.
    \end{equation}
\end{linenomath}

The strain energy function is used to calculate the second Piola-Kirchhoff stresses ${^t\mathbf{S}^{e}}$ and the stress-strain material property matrix $^t\mathbf{C}_\mathrm{NL}^{e}$ for the plane stress state via corresponding derivations and transformations. In contrast to other formulations for the compressible neo-Hookean material model, the required matrices for the plane stress state can be calculated directly for this formulation of the strain energy function, whereas otherwise an iterative calculation is necessary \cite{Campello2003}. For the detailed derivations of the second Piola-Kirchhoff stresses and the stress-strain material property matrix for the material model used here, see \cite{daSilva2020}.

Even with polyconvex material models, the convergence issues in elements with low stiffness can only be solved for relative small distortions of the CM. Additional stabilization is therefore needed. Many possible solutions are proposed in the literature. \textsc{Pedersen et al.} \cite{Pedersen2001} relax the convergence criterion for the Newton-Raphson iterations by deleting the DoFs of the nodes surrounded by elements with low stiffness from the convergence criterion. Another popular method is the additive hyperelasticity technique. For this, the finite elements are modelled with linear material behavior and additive hyperelastic material is temporarily applied to elements that become unstable \cite{Liu2017a, Luo2015} in order to stabilize them. In \cite{Kawamoto2009}, instead of the often used Newton-Raphson algorithm, the more robust Levenberg-Marquard algorithm is used for iterative solution in order to find a stable equilibrium. \textsc{Bruns et al.} \cite{Bruns2003} proposes an element removal and reintroduction algorithm. Here, elements of low stiffness are temporarily removed from the optimization in order to improve convergence. During optimization, these elements can also be added back into the design space. \textsc{Wang et al.} \cite{Wang2014c} introduces the energy interpolation scheme. Assuming that elements of low stiffness do not influence the structural behavior, they can be modelled arbitrarily. Therefore, they are modelled with linear FEM, which significantly improves convergence. 

We adapt the last mentioned approach for our optimization algorithm: elements with low stiffness ($\tilde{x}^{e}\approx 0$) are modelled with linear FEM, elements with high stiffness ($\tilde{x}^{e}=1$) with (geometrically and materially) nonlinear FEM. However, while the structure is forming, there are many elements with intermediate stiffnesses. An interpolation scheme is used to ensure that the transition in the modelling is smooth with increasing values of the design variables. Based on this method and according to \cite{XavierLeitao2019}, we model the parameterized internal forces and the parameterized tangent stiffness matrices as follows:

\begin{linenomath}
    \begin{equation}
        \label{Stiffness_interpolation}
        ^t\mathbf{K}^{e}(x^e)=x^{e}({^t\mathbf{K}^{e}_\mathrm{NL}}\upgamma^{e}+(1-\upgamma^{e})\mathbf{K}^{e}_\mathrm{L})
    \end{equation}
\end{linenomath}

\begin{linenomath}
    \begin{equation}
        \label{Load_interpolation}
        ^t\mathbf{f}^{e}(x^e)=x^{e}({^t\mathbf{f}^{e}_\mathrm{NL}}\upgamma^{e}+(1-\upgamma^{e})\mathbf{f}^{e}_\mathrm{L})
    \end{equation}
\end{linenomath}

They are therefore weighted with a weighting factor between linear and non-linear FEM. The weighting factor is calculated as in \cite{Wang2014c} with a Heaviside function:

\begin{linenomath}
    \begin{equation}
        \label{Heaviside_approximation_function}
        \upgamma^{e}=\frac{\tanh(\upbeta_{\upgamma}\upeta_{\upgamma})+\tanh(\upbeta_{\upgamma}(x^{e}-\upeta_{\upgamma}))}{\tanh(\upbeta_{\upgamma}\upeta_{\upgamma})+\tanh(\upbeta_{\upgamma}(1-\upeta_{\upgamma}))}
    \end{equation}
\end{linenomath}

When selecting the parameters in the Heaviside function, we follow the hints presented by \cite{Wang2014c} and select $\upbeta_{\upgamma}=500$ and $\upeta_{\upgamma}=0.01$.

Furthermore, the use of 4-node quadrilateral elements leads to numerical difficulties. Even with moderate deformations, the Newton-Raphson iterations do not converge if structures are only connected to each other via a node or an element. This must therefore be avoided during the complete iterative synthesis procedure. An effective method for this is the use of the density filter according to \cite{Bruns2001a,Bourdin2001}. The filtered design variables are calculated as follows:

\begin{linenomath}
    \begin{equation}
        \label{Density_filter}
        \tilde{x}^{e}=\frac{\sum_{i\in\upvartheta^{e}}w(\mathbf{P}^i)x^{i}}{\sum_{i\in\upvartheta^{e}}w(\mathbf{P}^i)}
    \end{equation}
\end{linenomath}

The centers $\mathbf{P}^i$ of the elements included in the calculation are located in a circular environment $\upvartheta_{e}$ in which the center $\mathbf{P}^e$ of the element under consideration $e$ lies in the center. The radius $R$ is specified. The linear weighting function used is:

\begin{linenomath}
    \begin{equation}
        \label{linear_weighting_function}
        w(\mathbf{P_i})=R-\Vert\mathbf{P}^i-\mathbf{P}^e\Vert
    \end{equation}
\end{linenomath}

This filter method also reduces the checkerboard patterns known from optimization with linear FEM.

Another well-known phenomenon is that the Newton-Raphson algorithm does not converge if the current displacement deviates too far from the displacement in the structural equilibrium \cite{Bathe2014}. For this purpose, external forces or displacement boundary conditions are usually applied incrementally and the Newton-Raphson iterations are performed for each increment. The topology optimization algorithm described in section \ref{Optimization_formulation} is also iterative. Some terminology should be noted: The Newton-Raphson iterations are not the iterations of the optimization algorithm. During an iteration of the optimization algorithm, several Newton-Raphson iterations are executed for each stationary point in order to find the structural equilibrium. The optimization algorithm described later has the property that the design variables and thus also the current displacement change only slightly between the iteration steps of the optimization algorithm. This property makes the optimization algorithm well suited for non-linear FEM, as the current displacements between the iteration steps of the optimization algorithm change only slightly. This means that in most iteration steps of the optimization algorithm, incremental application of the displacement can be omitted if the displacement in the structural equilibrium of the last iteration step is used as the initial value. The displacements only have to be applied incrementally for the Newton-Raphson iterations for each stationary point at the start of the optimization. For all further iterations of the optimization algorithm, the entire displacement is applied in one increment. However, convergence difficulties still occur in a few increments during optimization. In these iteration steps of the optimization algorithm, an automatic increase in the number of increments is used for the Newton-Raphson iterations. For convergence to a meaningful CM, it is also important to select a sufficiently small convergence tolerance for the convergence criterion.
\section{Design Problem}
\label{Design_Problem}

The tangent stiffness matrices in the stationary points $^t\mathbf{K},t=1...n$ represent, as already mentioned, a linearization over the specified deformation. All variables associated with the tangent stiffness matrices are labelled with $^{t}(*)$. For a path-generating CM with selective compliance, these tangent stiffness matrices must have properties that are explained below.

Since no external forces $^t\mathbf{r}$ may act on the passive DoFs, the tangent stiffness matrices can be divided into active and passive DoFs and then condensed to their active DoFs \cite{Gasch1989}:

\begin{linenomath}
    \begin{equation}
        \label{subdivision_stiffness_matrix}
        ^{t}\mathbf{K} {^{t}\mathbf{u}} = {^{t}\mathbf{r}} \iff
        \begin{bmatrix}
            ^{t}\mathbf{K}_{aa}&&^{t}\mathbf{K}_{ac}\\^{t}\mathbf{K}_{ca}&&^{t}\mathbf{K}_{cc}
        \end{bmatrix}
        \begin{bmatrix}
            ^{t}\mathbf{u}_{a}\\^{t}\mathbf{u}_{c}
        \end{bmatrix}
        =
        \begin{bmatrix}
            ^{t}\mathbf{r}_{a}\\\mathbf{0}
        \end{bmatrix}
    \end{equation}
\end{linenomath}
\begin{linenomath}
    \begin{equation}
        \label{condensation}
        (^{t}\mathbf{K}_{aa}-{^{t}\mathbf{K}_{ac}}~{^{t}\mathbf{K}_{cc}^{-1}}~{^{t}\mathbf{K}_{ca}}){^{t}\mathbf{u}_a}={^{t}\mathbf{\bar{K}}}~{^{t}\mathbf{u}_{a}}={^{t}\mathbf{r}_{a}}
    \end{equation}
\end{linenomath}

To determine the deformation behavior of a CM, the eigenbehavior has to be determined. The eigenmodes $^{t}\bar{\bm{\upchi}}_j,j=1...q$ and eigenvalues $^{t}\uplambda_j,j=1...q$ of each tangent stiffness matrix can be determined by solving the following eigenproblem:

\begin{linenomath}
    \begin{equation}
        \label{eigenvalue problem}
        ^{t}\mathbf{\bar{K}}~{^{t}\bar{\bm{\mathrm{\upchi}}}}={^{t}\uplambda}~{^{t}\bar{\bm{\mathrm{\upchi}}}}
    \end{equation}
\end{linenomath}

The eigenmodes have the following characteristics:

\begin{linenomath}
    \begin{equation}
        \label{normalization}
        ^{t}\bar{\bm{\upchi}}_j^\mathrm{T}~{^{t}\bar{\bm{\upchi}}_j}=1, \; j=1...q
    \end{equation}
\end{linenomath}

The eigenvalues are usually sorted in ascending order and the corresponding eigenmodes are arranged accordingly. In each stationary point, the first eigenmode is referred to as the kinematic eigenmode $^{t}\bar{\bm{\mathrm{X}}}_\mathrm{d}$ and all other $q-1$ eigenmodes belong to the subset of parasitic eigenmodes $^{t}\bar{\bm{\mathrm{X}}}_\mathrm{ud}$:

\begin{linenomath}
    \begin{equation}
        \label{subdivision_eigenmodes}
        ^{t}\bar{\bm{\mathrm{X}}}=
        \begin{bmatrix}
        ^{t}\bar{\bm{\mathrm{X}}}_\mathrm{d}&|&^{t}\bar{\bm{\mathrm{X}}}_\mathrm{ud}
        \end{bmatrix}, \;
        ^{t}\bar{\bm{\mathrm{X}}}_\mathrm{d}=
        \begin{bmatrix}
        ^{t}\bar{\bm{\upchi}}_{1}
        \end{bmatrix}, \;
        ^{t}\bar{\bm{\mathrm{X}}}_\mathrm{ud}=
        \begin{bmatrix}
        ^{t}\bar{\bm{\upchi}}_{2}&^{t}\bar{\bm{\upchi}}_{3}&\ldots&^{t}\bar{\bm{\upchi}}_{q}
        \end{bmatrix};
    \end{equation}
\end{linenomath}

Taking condition \eqref{normalization} into account, the following variables can be calculated in each stationary point of the CM using the first two eigenmodes:

\begin{linenomath}
    \begin{equation}
        \label{Primary_stiffness}
        ^{t}K_\mathrm{p}({^{t}\mathbf{\bar{K}}})={^{t}\uplambda_{1}}={^{t}\bar{\bm{\upchi}}_{1}^\mathrm{T}}~{^{t}\mathbf{\bar{K}}}~{^{t}\bar{\bm{\upchi}}_{1}}  
    \end{equation}
\end{linenomath}

\begin{linenomath}
    \begin{equation}
        \label{Secondary_stiffness}
        ^{t}K_\mathrm{s}({^{t}\mathbf{\bar{K}}})={^{t}\uplambda_{2}}={^{t}\bar{\bm{\upchi}}_{2}^\mathrm{T}}~{^{t}\mathbf{\bar{K}}}~{^{t}\bar{\bm{\upchi}}_{2}}    
    \end{equation}
\end{linenomath}

The primary stiffness $^{t}K_\mathrm{p}$ is calculated using the first eigenmode and corresponds to the first eigenvalue. The secondary stiffness $^{t}K_\mathrm{s}$ is calculated using the second eigenmode and corresponds to the second eigenvalue. These considerations are valid for CM with a pseudo-mobility of 1. The pseudo-mobility is a variable adapted for CM that is comparable to the mobility of conventional mechanisms \cite{Seltmann2022}. A pseudo-mobility of 1 means that only one kinematic eigenmode is present.

\begin{figure}[H]
\centering
\includegraphics[width=0.8\textwidth]{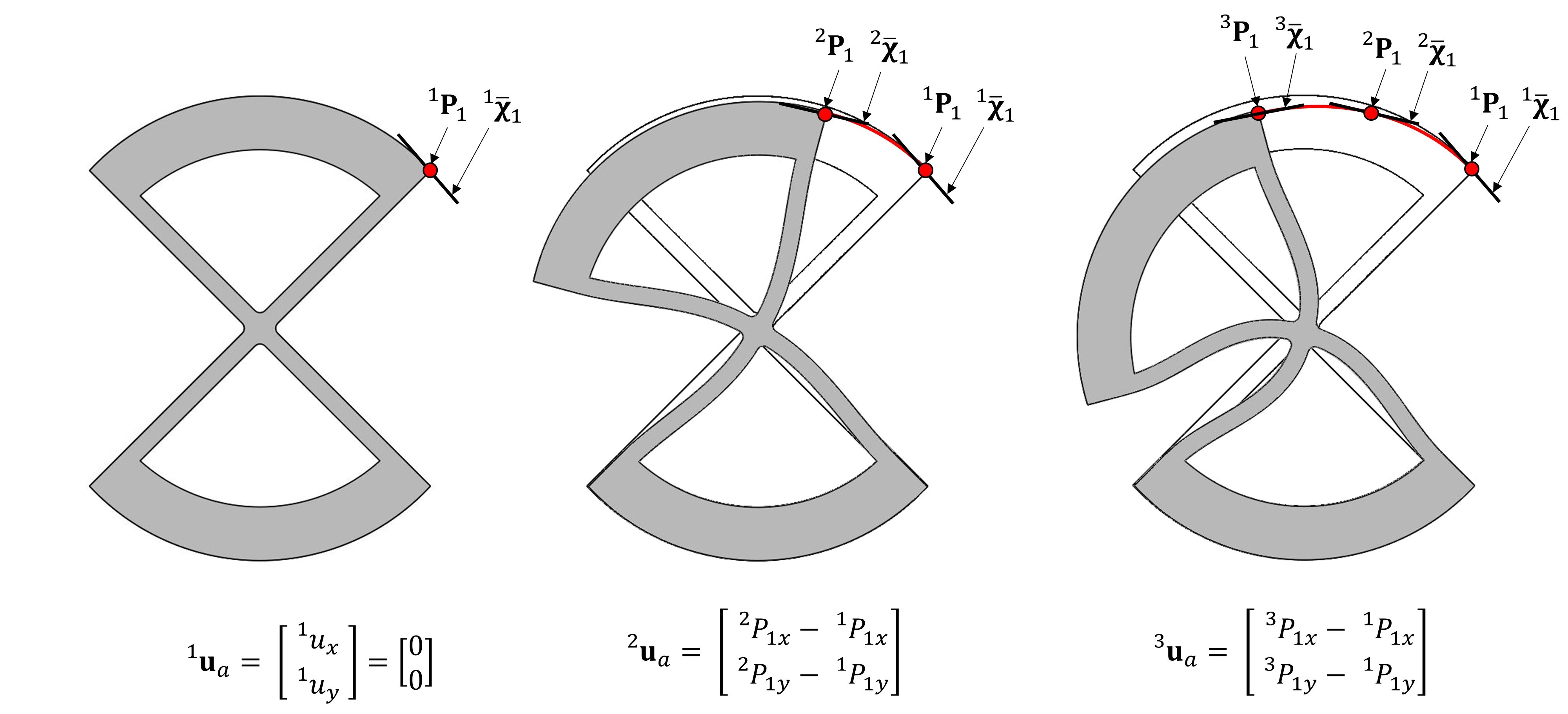}
\caption{Deformation path of a cartwheel hinge. \label{fig1}}
\end{figure}

The eigenmodes in the stationary points indicate all possible linearized deformations of the active DoFs, according to which the CM can be further deformed. The eigenvalues are a measure of the stiffness that is opposed to a further deformation of the CM deformed at the stationary point $t$ with $^t\mathbf{u}_a$ according to the respective associated eigenmode. The further deformation of the active DoFs will be a linear combination of the eigenmodes (kinematic and parasitic) for small additional deformations. If the ratio between secondary and primary stiffness (hereafter referred to as selectivity) becomes large, deformations corresponding to the parasitic eigenmodes become unfavorable, as a high stiffness must be overcome. The greater the selectivity becomes, the lower the proportion of parasitic eigenmodes in the resulting deformation of the active DoFs (the precision becomes higher). If the selectivity is infinitely large in all stationary points, the CM with selective compliance would exhibit the property of a conventional mechanism: its active DoFs would always deform along exactly the same deformation path regardless of the applied load. A sufficiently high selectivity can also be created for CMs (CMs with selective compliance). An example of this is shown in Figure~\ref{fig1}. Here, only the two DoFs in a single point is defined as active DoFs for clarity.

\unskip

For the optimization, it must also be required that the kinematic eigenmode approximates the desired tangent deformation mode at each stationary point $^t\mathbf{u}_a$. The required approximation and the high selectivity in all stationary points create a path that is energetically favorable for the deformation of the CM. Deviations from this deformation path are energetically unfavorable for the deformation. If the kinematic eigenmode in each stationary point is not identical to the desired tangent deformation mode after optimization, the minimum energy path will deviate from the desired kinematics. This deviation is referred to as accuracy and the generated minimum energy path as the natural kinematics of the CM \cite{Campanile2022}.

Stationary points that are elements of the desired kinematics are selected for the optimization. However, there are slight limitations to the practical feasibility of path-generating CM with selective compliance. The deformation behavior can only be optimized for a certain domain of definition that is covered by the stationary points. This is referred to below as the domain of definition of the CM. If forces act on the active DoFs in such a way that the deformation leaves this range, the desired behavior cannot be ensured. Therefore, the choice of the area of the natural kinematics to be optimized and thus the stationary points must be adapted to the application.

\section{Optimization formulation}
\label{Optimization_formulation}

\subsection{Problem statement}
\label{Problem_statement}

The optimization procedure developed for linear assumptions in \cite{Kirmse2021} serves as the basis for the procedure presented here. The optimization procedure must ensure that a high selectivity can be maintained at each stationary point. In addition, the first eigenmode in each stationary point must approximate to specified desired tangent deformation mode $^t\bar{\bm{\upvarphi}}$. By selecting the displacements in the stationary points $^t\mathbf{u}_a$ and the corresponding desired tangent deformation modes $^t\bar{\bm{\upvarphi}}$, the designer can specify the deformation paths of the active DoFs to be generated. A set of $n$ tangent stiffness matrices is calculated for the specified displacements. These are all parameterized to $^t\mathbf{K}(\mathbf{x})$ with the same set of design variables \eqref{Parametrization}  and can then be reduced to $^{t}\mathbf{\bar{K}}(\mathbf{x})$ for each stationary point using the equation \eqref{condensation}. The general optimization formulation, which summarizes all stationary points, is as follows

\begin{linenomath}
    \begin{equation}
        \label{Secondary_stiffness_general}
        \mathrm{max} \; {f{(\mathbf{x})}}=\sum_{t=1}^n {^t\upomega}K_\mathrm{s}({^t\mathbf{\bar{K}}(\mathbf{x})})
    \end{equation}
\end{linenomath}   
such that:
\begin{linenomath}
    \begin{equation}
        \label{Restriction_primary_stiffness_general1}
        \begin{matrix}
        \text{variant 1} \\ g(\mathbf{x})=\sum_{t=1}^n{^t\bar{\bm{\upvarphi}}^\mathrm{T}}~{^t\mathbf{\bar{K}}(\mathbf{x})}~{^t\bar{\bm{\upvarphi}}} -2 \upmu_\mathrm{g} \le 0
        \end{matrix}
    \end{equation}
\end{linenomath}      
\begin{align}
    \label{Restriction_primary_stiffness_general2}
    \begin{matrix}
        \text{variant 2: } \\ g_1(\mathbf{x})\ldots g_n(\mathbf{x})={^t\bar{\bm{\upvarphi}}^\mathrm{T}}~{^t\bar{\mathbf{K}}(\mathbf{x})}~{^t\bar{\bm{\upvarphi}}} -2~{^t\upmu} \le 0, t=1...n
    \end{matrix}
\end{align}
\begin{linenomath}
    \begin{equation}
        \label{Volume_general}
        m(\mathbf{x})=\sum_{e=1}^m (x^e)^{1/\upeta}-m V\le 0
    \end{equation}
\end{linenomath} 
\begin{linenomath}
    \begin{equation}
        \label{Restriction_x_general}
        x_\mathrm{l}^{e} \le x^{e} \le x_\mathrm{u}^{e}, e=1...m
    \end{equation}
\end{linenomath} 

The secondary stiffness is maximized using a sum function in each stationary point \eqref{Secondary_stiffness_general}. The individual secondary stiffnesses must be weighted accordingly so that no secondary stiffness dominates the optimization. The calculation of the weighting factors is described in section~\ref{Global_optimization}. The primary stiffnesses for all stationary points must be limited to a constant value. There are two variants for this. In variant~1, the primary stiffnesses for all stationary points are restricted as a sum with a value $\upmu_\mathrm{g}$. In variant~2, the primary stiffnesses are restricted individually with a set of constraints with the value $^t\upmu$. Further information on the choice of $\upmu$ can be found in \cite{Kirmse2021}. The volume is restricted in equation \eqref{Volume_general}. The calculation of the volume differs from the usual formulation in that it is not the sum of $x^e$, but the sum of $(x^{e})^{1/\upeta}$, where $\upeta$ represents a penalty factor. Due to this, the result tends towards 0-1 solutions (the design variables take on either minimum or maximum possible values). Further information on this special penalization procedure can be found in \cite{Seltmann2023} and \cite{Bruns2005}. The permissible values of the design variables are restricted in \eqref{Restriction_x_general} to avoid numerical problems.

The optimization problem presented is difficult to solve, so it is divided into two subproblems. For subproblem~1, an orthonormal base $^t\bar{\bm{\uppsi}}=[{^t\bar{\bm{\upvarphi}}},{^t\bar{\bm{\uppsi}}_1},{^t\bar{\bm{\uppsi}}_2},. ...,{^t\bar{\bm{\uppsi}}}_{q-1}]$ to the matrix $^{t}\mathbf{\bar{K}}(\mathbf{x}_{\mathrm{t0}})$, which is assumed to be constant, is computed using an appropriate optimization formulation. The orthonormal base $^t\bar{\bm{\uppsi}}$ represents an approximation to the eigenmodes of $^{t}\mathbf{\bar{K}}(\mathbf{x}_\mathrm{t0})$. In $^t\bar{\bm{\uppsi}}$, the first vector is the desired tangent deformation mode $^t\bar{\bm{\upvarphi}}$ for the respective stationary point. The remaining vectors are an approximation to the undesired deformation modes. The orthonormal base is then expanded to all structural DoFs to $^t\bm{\Uppsi}$. In subproblem~2, $^{t}\mathbf{K}(\mathbf{x})$ is varied for the expanded problem and the orthonormal base is seen as constant. The two subproblems are iterated consecutively. The optimized design variables $\mathbf{x}$ are used to calculate $\mathbf{x}_\mathrm{t0}$ for subproblem~1. An iteration step includes a run of subproblem~1 and 2 and is denoted by $s$.

The presented optimization problem is reduced to a topology optimization for linear CM for one stationary point with $^1\mathbf{u}_a=\mathbf{0}$.

\subsection{Subproblem~1: calculation of the orthonormal base and expansion}

The following optimization problem must be solved to calculate the orthonormal base for each stationary point:

\begin{align}
    \label{Golub}
    \mathrm{min} \; {^tf({^t\bar{\bm{\uppsi}}_j})}= {^t\bar{\bm{\uppsi}}_j^\mathrm{T}}~ {^t\mathbf{\bar{K}}(\mathbf{x}_\mathrm{t0})}~{^t\bar{\bm{\uppsi}}_j}
\end{align}
such that:
\begin{align}
    \label{Golub_Con1}
    \begin{aligned}
      ^tg_1({^t\bar{\bm{\uppsi}}_j}) &= {^t\bar{\bm{\upvarphi}}^\mathrm{T}}~{^t\mathbf{\bar{K}}(\mathbf{x}_\mathrm{t0})}~{^t\bar{\bm{\uppsi}}_j} =0
    \end{aligned} 
\end{align}
\begin{align*}
    \label{Golub_Con2}    
    ^tg_{2}({^t\bar{\bm{\uppsi}}_j}) &= {^t\bar{\bm{\uppsi}}_1^\mathrm{T}}~{^t\mathbf{\bar{K}}(\mathbf{x}_\mathrm{t0})}~{^t\bar{\bm{\uppsi}}_j =0} ,\;j=2 \\
\end{align*}
\begin{align}
    \left. 
    \begin{aligned}
        ^tg_{2}({^t\bar{\bm{\uppsi}}_j}) &= {^t\bar{\bm{\uppsi}}_1^\mathrm{T}}~ {^t\mathbf{\bar{K}}(\mathbf{x}_\mathrm{t0})}~{^t\bar{\bm{\uppsi}}_j} =0 \\ &\vdots \\ ^tg_{j}({^t\bar{\bm{\uppsi}}_j}) &= {^t\bar{\bm{\uppsi}}_{j-1}^T}~{^t\mathbf{\bar{K}}(\mathbf{x}_\mathrm{t0})}~{^t\bar{\bm{\uppsi}}_j} =0 \\
    \end{aligned} 
    \right\} j>2
\end{align}
\begin{align}
    \label{Golub_Con3}
    ^th({^t\bar{\bm{\uppsi}}_j}) = {^t\bar{\bm{\uppsi}}_j^\mathrm{T}}~{^t\bar{\bm{\uppsi}}_j} =1
\end{align}

This optimization problem is taken from \cite{Kirmse2021} and must be solved recursively for $j=1...q-1$. A computationally efficient solution to this problem is described in \cite{Golub1973b}. All vectors $^t\bar{\bm{\uppsi}}_j$ of the orthonormal base are calculated with the help of a substitute eigenvalue problem. To position $^t\bar{\bm{\uppsi}}_j$ correctly in $^t\bar{\bm{\Uppsi}}$, they are ordered in ascending order by their function value in \eqref{Golub}.

Then all vectors of the base $^t\bar{\bm{\Uppsi}}$ are expanded to all structural DoFs:

\begin{align}
    \label{nr18}
    ^t\bm{\upvarphi}=
    \begin{bmatrix}
        {^t\bar{\bm{\upvarphi}}} \\ -{^t\mathbf{K}_{aa}^{-1}}~{^t\mathbf{K}_{ca}}~{^t\bar{\bm{\upvarphi}}} 
    \end{bmatrix}
    ; \;
    {^t\bm{\uppsi}_j}=
    \begin{bmatrix}
        {^t\bar{\bm{\uppsi}}_j} \\ -{^t\mathbf{K}_{aa}^{-1}}~{^t\mathbf{K}_{ca}}~{^t\bar{\bm{\uppsi}}_j} 
    \end{bmatrix},\; j=1 \ldots q-1
\end{align}

\subsection{Subproblem~2: updating the design variables}

Subproblem~2 corresponds in general to the optimization problem \eqref{Secondary_stiffness_general} to \eqref{Restriction_x_general} rewritten for the uncondensed system:

\begin{align}
    \label{Maximization_secondary_stiffness}
    \begin{matrix}
    \mathrm{max} \; f(\mathbf{x})=\sum_{t=1}^n {^t\upomega}~{^t\bm{\uppsi}_1^\mathrm{T}}~{^t\mathbf{K}(\mathbf{x})}~{^t\bm{\uppsi}_1}
    \end{matrix}
\end{align}

such that:

\begin{align}
    \label{Restriction_primary_stiffness1}
    \begin{matrix}
          \text{variant~1: } \\ g_1(\mathbf{x})=\sum_{t=1}^n {^t\bm{\upvarphi}^\mathrm{T}}~{^t\mathbf{K}(\mathbf{x})}~{^t\bm{\upvarphi}} -2\upmu_\mathrm{G} \le 0
    \end{matrix}
\end{align}
\begin{align}
    \label{Restriction_primary_stiffness2}
    \begin{matrix}
        \text{variant~2: } \\ g_1(\mathbf{x})\ldots g_n(\mathbf{x})={^t\bm{\upvarphi}^\mathrm{T}}~{^t\mathbf{K}(\mathbf{x})}~{^t\bm{\upvarphi}} -2~{^t\upmu} \le 0, t=1...n
    \end{matrix}
\end{align}
\begin{align}
    \label{Restriction_mode_swap}
    \begin{matrix}
    k_1(\mathbf{x})\ldots k_{n\cdot (l-1)}(\mathbf{x})={^t\bm{\uppsi}_1^\mathrm{T}}~{^t\mathbf{K}(\mathbf{x})}~{^t\bm{\uppsi}_1}-{^t\bm{\uppsi}_j^\mathrm{T}}~{^t\mathbf{K}(\mathbf{x})}~{^t\bm{\uppsi}_j}\le 0, j=2...l, t=1...n
    \end{matrix}
\end{align}
\begin{align}
    \label{Volume_restriction}
    \begin{matrix}
    m(\mathbf{x})=\sum_{e=1}^m \frac{(x^e_\mathrm{t0})^{1/\upeta}}{x^e_\mathrm{t0}}x^e  - mV \le 0
    \end{matrix}
\end{align}
\begin{align}
    \label{Restriction_boundaries_global}
    x_\mathrm{l}^{e} \le x^{e} \le x_\mathrm{u}^{e}, \; e=1...m
\end{align}

The equation \eqref{Maximization_secondary_stiffness} corresponds to the equation \eqref{Secondary_stiffness_general}. The secondary stiffness is maximized simultaneously for all stationary points using the second vector $^t\bar{\bm{\uppsi}}_1$ of the respective orthonormal bases. This vector approximates the first undesired deformation mode of the stiffness matrix. The calculation of the associated weighting factors in the sum function is described in the following section. The constraint of the primary stiffness for variant~1 and variant~2 in the equations \eqref{Restriction_primary_stiffness1} and \eqref{Restriction_primary_stiffness2} correspond to the equations \eqref{Restriction_primary_stiffness_general1} and \eqref{Restriction_primary_stiffness_general2} for the uncondensed system. The first undesired deformation mode should not swap its rank with another undesired deformation mode at any stationary point during the optimization, as this could destabilize the optimization algorithm. To prevent this, a corresponding constraint is defined in \eqref{Restriction_mode_swap} for each stationary point for a user-defined number $l$ of undesired deformation modes \cite{Kirmse2021}. The volume constraint in equation \eqref{Volume_restriction} corresponds to equation \eqref{Volume_general}. This differs from equation \eqref{Volume_general} in that the design variables $\mathbf{x}_\mathrm{t0}$, for which the orthnonormal base was calculated, are included in the equation for penalization. This method was described in \cite{Seltmann2023} and is permissible as long as the change in the design variables between two successive iteration steps is kept as small as possible. Due to this linearization, the optimization problem \eqref{Maximization_secondary_stiffness} to \eqref{Restriction_boundaries_global} can be solved easily. Here, the simplex algorithm according to \cite{Nelder1965} is used for the solution.

\subsection{Global optimization procedure}
\label{Global_optimization}

As already mentioned, an iteration step $s$ of the optimization procedure includes the determination of the tangent stiffness matrices for all stationary points. The objective function and constraints (subproblems~1 and 2) are subsequently set up for each stationary point. The design variables are optimized as a next step. The design variables $\mathbf{x}$ are then filtered using equation \eqref{Density_filter} to $\tilde{\mathbf{x}}$.

In order to converge to a solution, it is necessary to keep the change in the design variable between the iteration steps as small as possible. To do this, the start value for the next iteration step $s+1$ is calculated as a linear combination of the start value of the current iteration step $\mathbf{x}_{t0}(s)$ and the solution of the current iteration step $\tilde{\mathbf{x}}(s)$:

\begin{linenomath}
    \begin{equation}
        \label{Update_design_variable}
        \mathbf{x}_\mathrm{t0}(s+1)=\upkappa\mathbf{x}_\mathrm{t0}(s)+(1-\upkappa)\tilde{\mathbf{x}}(s) 
    \end{equation}
\end{linenomath}

The optimization problem \eqref{Maximization_secondary_stiffness} to \eqref{Restriction_boundaries_global} generates CM with concentrated compliance for both linear and nonlinear FEM. This is inherent to the optimization formulation and is explained in \cite{Seltmann2023}. Therefore, it is necessary to modify it accordingly to obtain CM with selective compliance. In \cite{Seltmann2023} an adaptive volume constraint is used to solve this problem. In the presence of a constraint for the primary stiffness, the adaptive volume constraint ensures that a CM with selective compliance is obtained. This methodology eliminates the need to use a weighting factor in the objective function, for example. However, problems can arise with complex design examples. Therefore, a different method is used here, but it is based on the same idea: First, for a certain number $I_\mathrm{T}$ of iteration steps $s$, the optimization is performed according to the optimization problem in subproblem~2. Then the volume constraint \eqref{Volume_restriction} is deleted and the objective function in equation \eqref{Maximization_secondary_stiffness} is replaced by

\begin{linenomath}
    \begin{equation}
        \begin{split}
            \label{Max_secondary_stiffness_min_volume}
            \mathrm{max} \; {f{(\mathbf{x})}}&=f_1(\mathbf{x})-\upzeta\upomega_\mathrm{V}(s) f_2(\mathbf{x})\\&=\sum_{t=1}^n {^t\upomega(s)}~{^t\bm{\uppsi}^\mathrm{T}}~{^t\mathbf{K}(\mathbf{x})}~{^t\bm{\uppsi}}-\upzeta\upomega_\mathrm{V}(s)\sum_{e=1}^m \frac{(x^e_\mathrm{t0}(s))^{1/\upeta}}{x^e_\mathrm{t0}(s)}x^e(s),s>It
        \end{split}
    \end{equation}
\end{linenomath} 

This minimizes the volume in addition to maximizing the secondary stiffness. To ensure that one requirement does not dominate the other, it is necessary to introduce an additional weighting factor $\upomega_\mathrm{V}(s)$:

\begin{linenomath}
    \begin{equation}
        \label{Weight factor_volume}
        \upomega_\mathrm{V}(s) = \frac{-f_2(\mathbf{x}_\mathrm{t0}(s))}{f_1(\mathbf{x}_\mathrm{t0}(s))}
    \end{equation}
\end{linenomath}   

This is updated for each iteration step, resulting in an approximately constant weighting of the two requirements across all iteration steps. The factor $\upzeta<1$ is also introduced. This can be used to adjust the weighting of the volume in relation to the secondary stiffness. If the volume is weighted too high, the limit for the primary stiffness cannot be reached. If the weighting is too low, CMs with concentrated compliance are created. This version requires the selection of a weighting factor $\upzeta$ for each design problem. The weighting factors $^t\upomega(s)$ are also calculated for each iteration step so that the secondary stiffnesses are weighted equally in each stationary point:

\begin{linenomath}
    \begin{equation}
        \label{Weight factor_secondary stiffness}
        ^t\upomega(s) = \frac{({^t\bm{\uppsi}^\mathrm{T}}~{^t\mathbf{K}(\mathbf{x}_\mathrm{t0}(s))}~{^t\bm{\uppsi}})^{-1}}{\sum_{t1=1}^{n1} (({^{t1}\bm{\uppsi}^\mathrm{T}}~{^{t1}\mathbf{K}(\mathbf{x}_\mathrm{t0}(s))}~{^{t1}\bm{\uppsi}})^{-1})},t=1...n
    \end{equation}
\end{linenomath}   

Other weightings are also possible, but these are not considered further. The iterations are run through until the following convergence criterion is satisfied: If the sum of the design variables

\begin{linenomath}
    \begin{equation}
        \label{Volume}
        V_\mathrm{s}(s)=\sum_{e=1}^m x^e(s) 
    \end{equation}
\end{linenomath}

remains within a certain range $\upepsilon$ for a certain number of iteration steps $n$, the optimization procedure is terminated:

\begin{linenomath}
    \begin{equation}
        \label{Convergence_criterion}
        (1-\updelta_\mathrm{C}) V_\mathrm{s}(s-s_\mathrm{C}-1)\leq\{V_\mathrm{s}(i)\mid i=s-s_\mathrm{C} \ldots s\}\leq (1+\updelta_\mathrm{C}) V_\mathrm{s}(s-s_\mathrm{C}-1)
    \end{equation}
\end{linenomath}

\begin{figure}[H]
    \centering
    \includegraphics[width=0.8\textwidth]{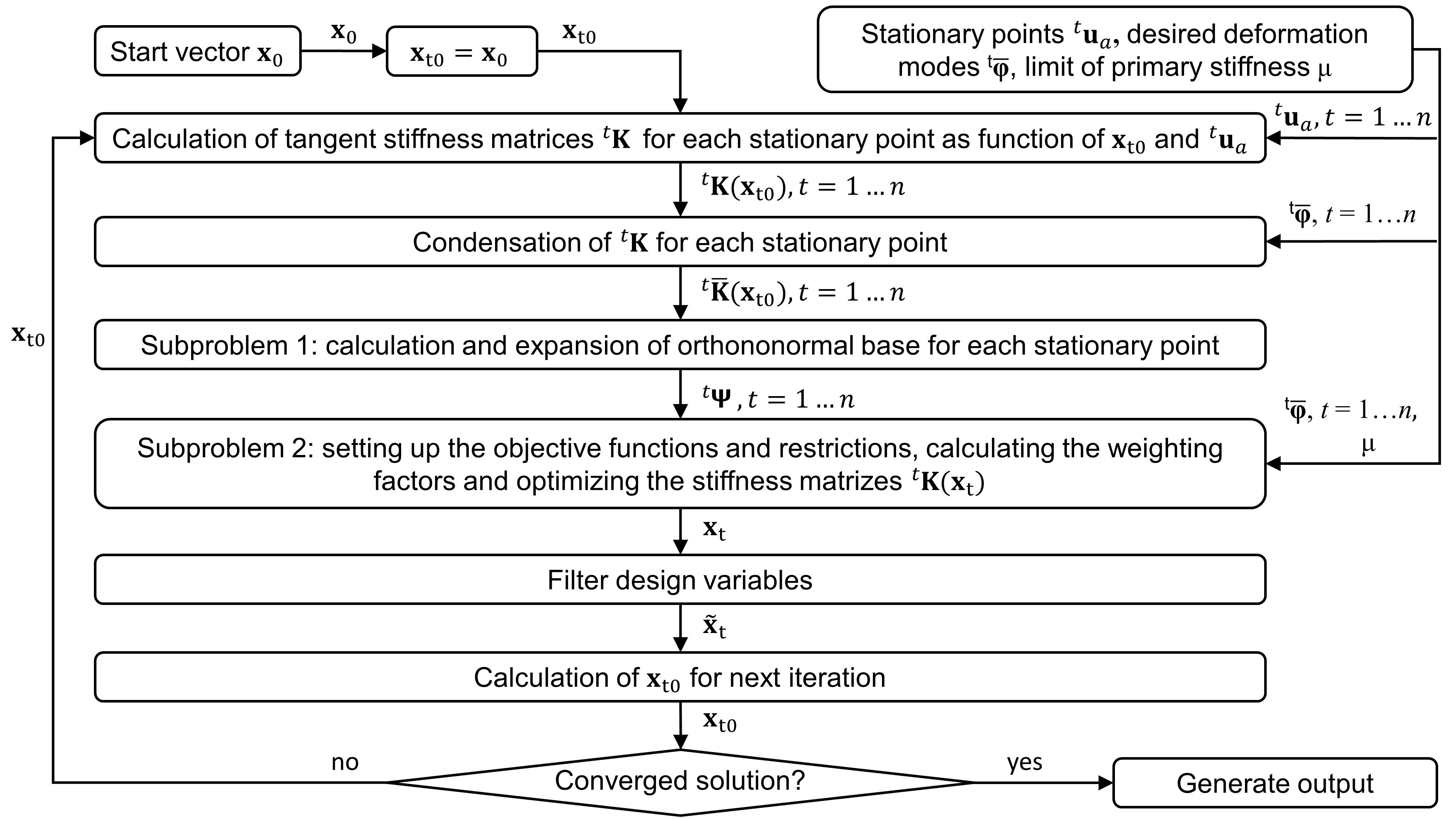}
    \caption{Global iteration procedure\label{fig2}}
    \end{figure}

The entire optimization procedure is shown in Figure \ref{fig2}.

    \unskip
\section{Design examples}

The optimization algorithm presented was implemented in MATLAB and tested on three different design examples for both linear and path-generating CM in order to perform a comparison. First, an inverter mechanism was implemented, which is a frequently used example in the literature. This inverts the direction of the input DoF displacement. The ratio of input to output DoF displacement was specified for the path-generating inverter for certain stationary points. A pivot joint was used as a second design example. This example can be used to show that it is possible to design CMs whose active DoFs deform according to predefined deformation paths in the plane. The third example is a shape-adaptive structure. The free surface of this structure should deform with increasing deformation in the form of a sinusoidal curve whose amplitude is increased. Finally, it can be shown that the presented optimization algorithm can also be used to calculate shape-adaptive path-generating structures, which to the best of the authors' knowledge has not yet been demonstrated in the literature.

A suitable design space must be defined for all design examples before optimization. This is filled with bilinear quadrilateral elements. In all design examples, the elements have a thickness of one. The dimensions of the side lengths differ in the various design examples. The same material properties are chosen for all design examples: Nylon with a modulus of elasticity of $E_\mathrm{M} = 3$ GPa and a Poissons's ratio of $\upnu_\mathrm{M} = 0.4$ is chosen as the material. Furthermore, the same values are applied to all entries of \eqref{Restriction_boundaries_global}:

\begin{linenomath}
    \begin{equation}
        \label{lower_limits}
        x_\mathrm{l}^{e}=x_\mathrm{l}, \; e=1...m
    \end{equation}
\end{linenomath}
\begin{linenomath}
    \begin{equation}
        \label{upper_limits}
        x_\mathrm{u}^{e}=x_\mathrm{u}, \; e=1...m
    \end{equation}
\end{linenomath}

All parameters that are selected in the same way for the optimizations for linear and path-generating CM for all design examples are listed in Table~\ref{tab1}.

\begin{table}[H] 
    \caption{Optimization parameters chosen to be the same for all design examples\label{tab1}}
    \newcolumntype{C}{>{\centering\arraybackslash}X}
    \begin{tabularx}{\textwidth}{CCC}
    \toprule
    \textbf{Optimization parameter}	& \textbf{Selected value}	& \textbf{In equation}\\
    \midrule
    $\upepsilon_\mathrm{F}$		& $1 \times 10^{-6}$			& \eqref{Residuum}\\
    $E_\mathrm{M}$		& $3$ GPa			& \eqref{Lame_parameters} \\
    $\upnu_\mathrm{M}$		& $0.4$			& \eqref{Lame_parameters}\\
    $\upbeta_\upgamma$		& $500$			& \eqref{Heaviside_approximation_function}\\
    $\upeta_\upgamma$		& $0.01$			& \eqref{Heaviside_approximation_function}\\
    $R$		& $1.5$			& \eqref{linear_weighting_function}\\
    $\upeta$		& $3$			& \eqref{Volume_restriction}\\
    $V$		& $0.3$			& \eqref{Volume_restriction}\\
    $\upkappa$		& $0.99$			& \eqref{Update_design_variable}\\
    $I_\mathrm{T}$		& $1000$			& \eqref{Max_secondary_stiffness_min_volume}\\
    $\updelta_C$		& $0.001$			& \eqref{Convergence_criterion}\\
    $s_\mathrm{C}$		& $500$			& \eqref{Convergence_criterion}\\
    $x_\mathrm{l}$		& $1 \times 10^{-9}$			& \eqref{lower_limits}\\
    $x_\mathrm{u}$		& $1$			& \eqref{upper_limits}\\
    \bottomrule
    \end{tabularx}
    \end{table}
    \unskip

A value of $0.5$ is selected as the starting value for all optimizations for all design variables $x^e$ and 10 increments for the Newton-Raphson iterations are selected for the first iteration for all design examples. 

\subsection{Inverter}

\begin{figure}[H]
    \centering
    \includegraphics[width=0.8\textwidth]{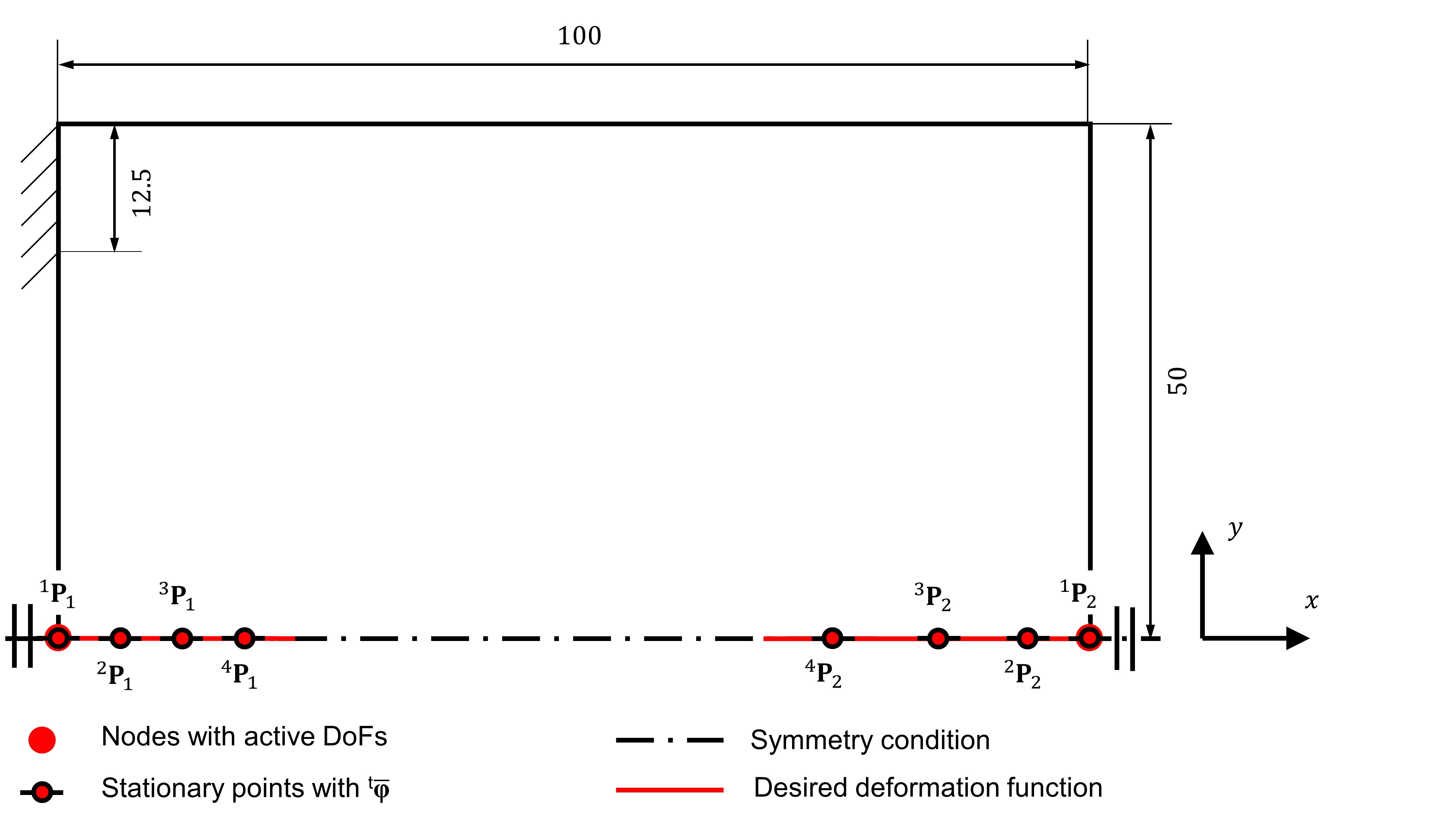}
    \caption{Inverter: design space with stationary points, desired deformation paths and the highlighted points with active DoFs \label{fig3}}
    \end{figure}

For the inverter, the same design space, the same clamping and the same number of elements are used for the linear and the path-generating CM. The design space is filled with 200 elements in the x-direction and 100 elements in the y-direction. The clamping, the symmetry condition used and the dimensions of the design space are shown in Figure~\ref{fig3}.

    \unskip

The structure has two active DoFs. These are the x-displacements of the two points highlighted in Figure~\ref{fig3}. The $x$-displacement at point $^tP_1$ is referred to as $u_{a1}$ and the $x$-displacement at point $^tP_2$ as $u_{a2}$. For better convergence of the Newton-Raphson iterations, the horizontally neighboring nodes are used instead of the nodes directly at the edge of the design space. The two active DoFs should perform an opposite displacement. Two cases with different desired deformation functions are optimized for the path-generating CMs. These should lead to the same stationary point and desired tangent deformation mode for $\upalpha=0$, which is used for the optimization of the linear CM. In the first case (inverter~a), the deformation function

\begin{equation}
    \label{desired_motion_nonlinear_Inverter1}
    \begin{bmatrix}
        \bar{\upxi}_1\\\bar{\upxi}_2
    \end{bmatrix}
    =
    \begin{bmatrix}
        \upalpha\\-\upalpha
    \end{bmatrix}
\end{equation}

is defined. In the second case (inverter~b), the deformation function is

\begin{equation}
    \label{desired_motion_nonlinear_Inverter2}
    \begin{bmatrix}
        \bar{\upxi}_1\\\bar{\upxi}_2
    \end{bmatrix}
    =
    \begin{bmatrix}
        \upalpha\\-\frac{1}{4}\upalpha^2-\upalpha
    \end{bmatrix}
\end{equation}

\begin{figure}[H]
    \centering
    \includegraphics[width=0.8\textwidth]{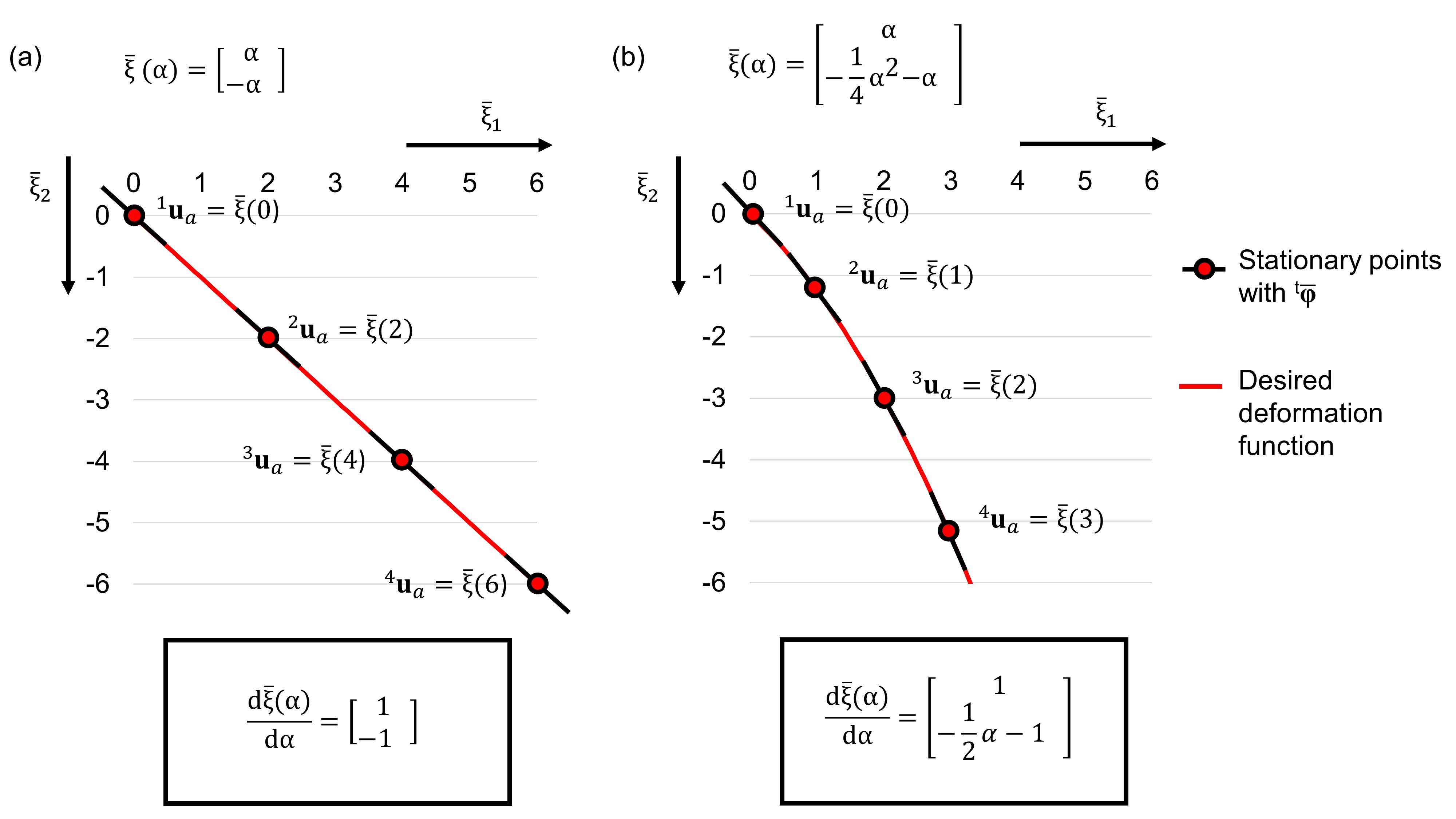}
    \caption{Inverter: functions describing the displacement ratio between the two active DoFs and their derivatives \label{fig4}}
    \end{figure}

The deformation functions and the associated derivatives, as well as the stationary points selected for the optimization for inverter~a and inverter~b are shown in Figure~\ref{fig4}.

    \unskip

The stationary points $^{t}\mathbf{u}_a$ can be determined directly from the functions. The corresponding desired tangent deformation modes $^{t}\bar{\bm{\upvarphi}}$ in the stationary points are determined by calculating the normalized derivatives. For inverter~a, the stationary points are thus defined as follows: 

\begin{linenomath}
    \begin{equation}
        \begin{split}
            \label{ua1a}
            ^{1}\mathbf{u}_a=
            \begin{bmatrix}
                0 \\ 0
            \end{bmatrix},\;
            ^{2}\mathbf{u}_a=
            \begin{bmatrix}
                2 \\ -2
            \end{bmatrix},\;
            ^{3}\mathbf{u}_a=
            \begin{bmatrix}
                4 \\ -4
            \end{bmatrix},\;
            ^{4}\mathbf{u}_a=
            \begin{bmatrix}
                6 \\ -6
            \end{bmatrix}\;
        \end{split}
    \end{equation}
\end{linenomath} 

In this example, the desired tangent deformation modes are the same at all stationary points: 

\begin{linenomath}
    \begin{equation}
        \begin{split}
            \label{phia}
            ^{1-4}\bar{\bm{\upvarphi}}&=
            \begin{bmatrix}
                0.7071 \\ -0.7071
            \end{bmatrix}\;
        \end{split}
    \end{equation}
\end{linenomath} 

For inverter~b the following stationary points and desired tangent deformation modes result:

\begin{linenomath}
    \begin{equation}
        \begin{split}
            \label{ua1b}
            ^{1}\mathbf{u}_a=
            \begin{bmatrix}
                0 \\ 0
            \end{bmatrix},\;
            ^{2}\mathbf{u}_a=
            \begin{bmatrix}
                1 \\ -5/4
            \end{bmatrix},\;
            ^{3}\mathbf{u}_a=
            \begin{bmatrix}
                2 \\ -3
            \end{bmatrix},\;
            ^{4}\mathbf{u}_a=
            \begin{bmatrix}
                3 \\ -21/4
            \end{bmatrix}\;
        \end{split}
    \end{equation}
\end{linenomath}

\begin{linenomath}
    \begin{equation}
        \begin{split}
            \label{phi1b}
            ^{1}\bar{\bm{\upvarphi}}=
            \begin{bmatrix}
                0.7071 \\ -0.7071
            \end{bmatrix},\;
            ^{2}\bar{\bm{\upvarphi}}=
            \begin{bmatrix}
                0.5547 \\ -0.8321
            \end{bmatrix},\;
            ^{3}\bar{\bm{\upvarphi}}=
            \begin{bmatrix}
                0.4472 \\ -0.8944
            \end{bmatrix},\;
            ^{4}\bar{\bm{\upvarphi}}=
            \begin{bmatrix}
                0.3714 \\ -0.9285
            \end{bmatrix}\;
        \end{split}
    \end{equation}
\end{linenomath} 

Further specific optimization parameters must be defined for the inverters, which are not defined in Table~\ref{tab1}. The path-generating inverters are optimized for variant~1 (see \eqref{Restriction_primary_stiffness1}) and variant~2 (see \eqref{Restriction_primary_stiffness2}). The associated limits for the primary stiffness are set to $\upmu_\mathrm{g}=0.5$ for the linear and the path-generating inverter with variant~2 and to $\upmu_\mathrm{g}=2$ for the path-generating inverter with variant~1 as a sum for four stationary points. As the inverters only have two active DoFs, there are only two possible eigenmodes. Therefore, $l$ is set to one in \eqref{Restriction_mode_swap}. The weighting for the volume constraint is set to $\upzeta=0.7$ in the objective function \ref{Max_secondary_stiffness_min_volume} for all inverters.

\subsection{Pivot joint}

\begin{figure}[H]
    \centering
    \includegraphics[width=0.8\textwidth]{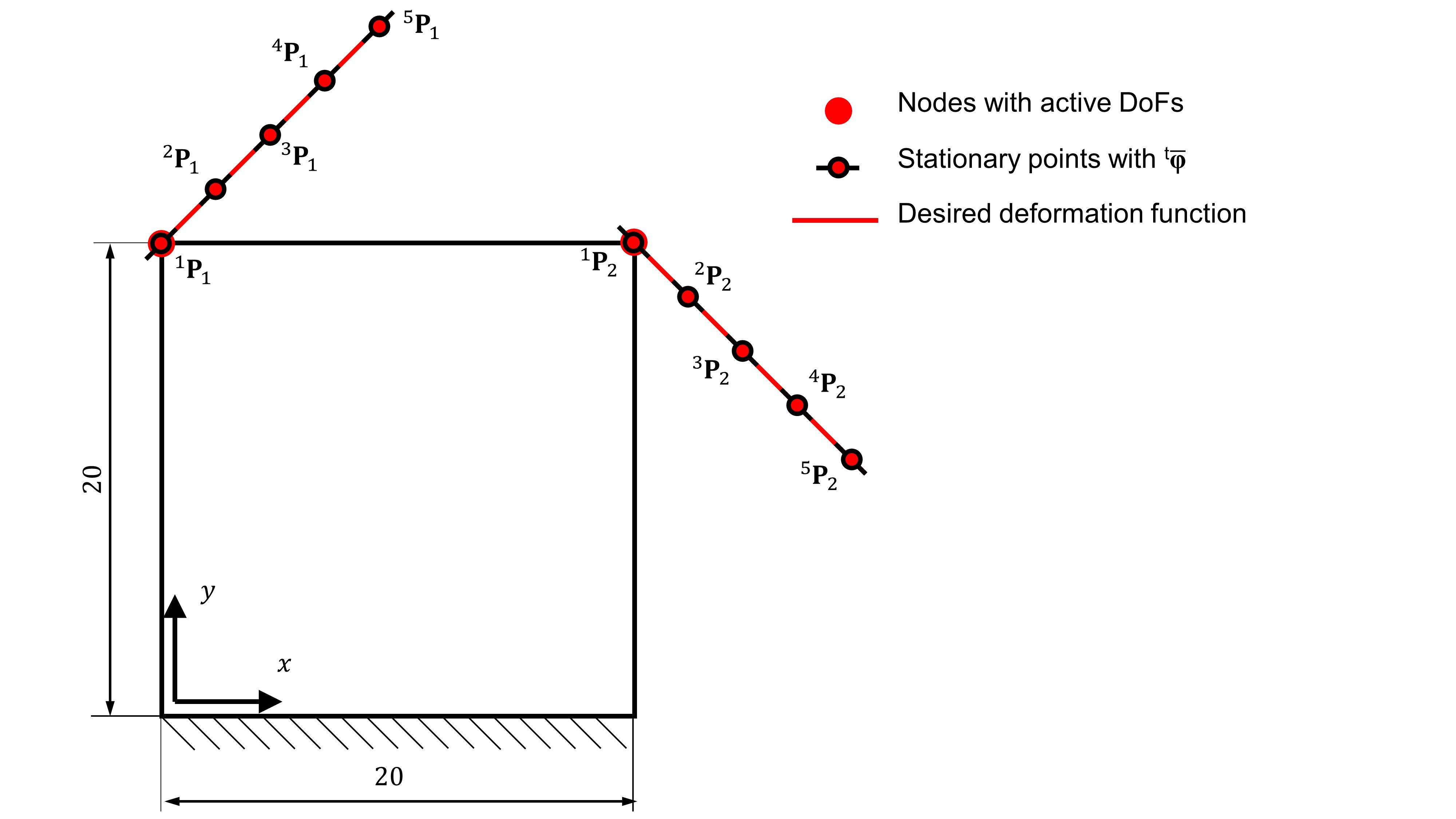}
    \caption{Pivot joint: design space with stationary points, desired deformation paths and the highlighted points with active DoFs \label{fig5}}
    \end{figure}

The defined design space and the boundary conditions for the linear and path-generating pivot joint are shown in Figure~\ref{fig5}. The design space is filled with 150 elements in both the x- and y-directions. Four active DoFs are defined in the two points highlighted in \ref{fig5}. The $x$-displacement at point $^tP_1$ is defined as $u_{a1}$, the $y$-displacement in point $^tP_1$ as $u_{a2}$, the $x$-displacement at point $^tP_2$ as $u_{a3}$ and the $y$-displacement at point $^tP_2$ as $u_{a4}$. Again, instead of using the outer corner nodes of the design space, we use the nearest diagonal node in the inner part of the design space. The desired tangent deformation modes are defined in the x- and y-directions of the two highlighted points. Point $\mathbf{P}_1$ should move on line~1, which is angled at 45°. Point $\mathbf{P}_2$ should move on line~2, which is angled at -45°, with the same amount as $\mathbf{P}_1$ on line~1. With these requirements, the desired deformation function can be created:

\begin{equation}
    \label{desired_motion_nonlinear_Pivot}
    \begin{bmatrix}
        \bar{\upxi}_1\\\bar{\upxi}_2\\\bar{\upxi}_3\\\bar{\upxi}_4
    \end{bmatrix}
    =
    \begin{bmatrix}
        \upalpha\\\upalpha\\\upalpha\\-\upalpha
    \end{bmatrix}
\end{equation}

Using this, the values for $^{t}\mathbf{u}_a$ and $^{t}\bar{\bm{\upvarphi}}$ can be determined again. For the choice of $\upalpha$, it is assumed that the stationary points are evenly distributed on line~1 (see figure~\ref{fig5}):

\begin{linenomath}
    \begin{equation}
        \begin{split}
            \label{ua2}
            ^{1}\mathbf{u}_a&=
            \begin{bmatrix}
                0 \\ 0 \\ 0 \\ 0
            \end{bmatrix},\;
            ^{2}\mathbf{u}_a=
            \begin{bmatrix}
                0.25 \\ 0.25 \\ 0.25 \\ -0.25
            \end{bmatrix},\;
            ^{3}\mathbf{u}_a=
            \begin{bmatrix}
                0.5 \\ 0.5 \\ 0.5 \\ -0.5
            \end{bmatrix},\;
            ^{4}\mathbf{u}_a=
            \begin{bmatrix}
                0.75 \\ 0.75 \\ 0.75 \\ -0.75
            \end{bmatrix},\;
            ^{5}\mathbf{u}_a=
            \begin{bmatrix}
                1 \\ 1 \\ 1 \\ -1
            \end{bmatrix}
        \end{split}
    \end{equation}
\end{linenomath} 

\begin{linenomath}
    \begin{equation}
        \begin{split}
            \label{phi2}
            ^{1-5}\bar{\bm{\upvarphi}}&=
            \begin{bmatrix}
                0.5 \\ 0.5 \\ 0.5 \\ -0.5
            \end{bmatrix}
        \end{split}
    \end{equation}
\end{linenomath} 

    \unskip

In order to plot the desired deformation function the graph is displayed separately for each individual point $P_i$ with active DoFs. The optimization for the linear CM can be carried out if only stationary point~1 is included in the optimization. The optimization for the path-generating CM is carried out for variant~1. The limit for the total primary stiffness (see \ref{Restriction_primary_stiffness1}) is set to $\upmu_\mathrm{g}=0.4$ for the linear CM and to $\upmu_\mathrm{g}=2$ for the path-generating CM as a sum for five stationary points. The parameter $l$ in \eqref{Restriction_mode_swap} is set to three. The weighting for the volume constraint is set to $\upzeta=0.7$ in the objective function \ref{Max_secondary_stiffness_min_volume} for the linear CM and to $\upzeta=0.3$ for the path-generating CM.

\subsection{Shape-adaptive structure}

\begin{figure}[H]
    \centering
    \includegraphics[width=0.8\textwidth]{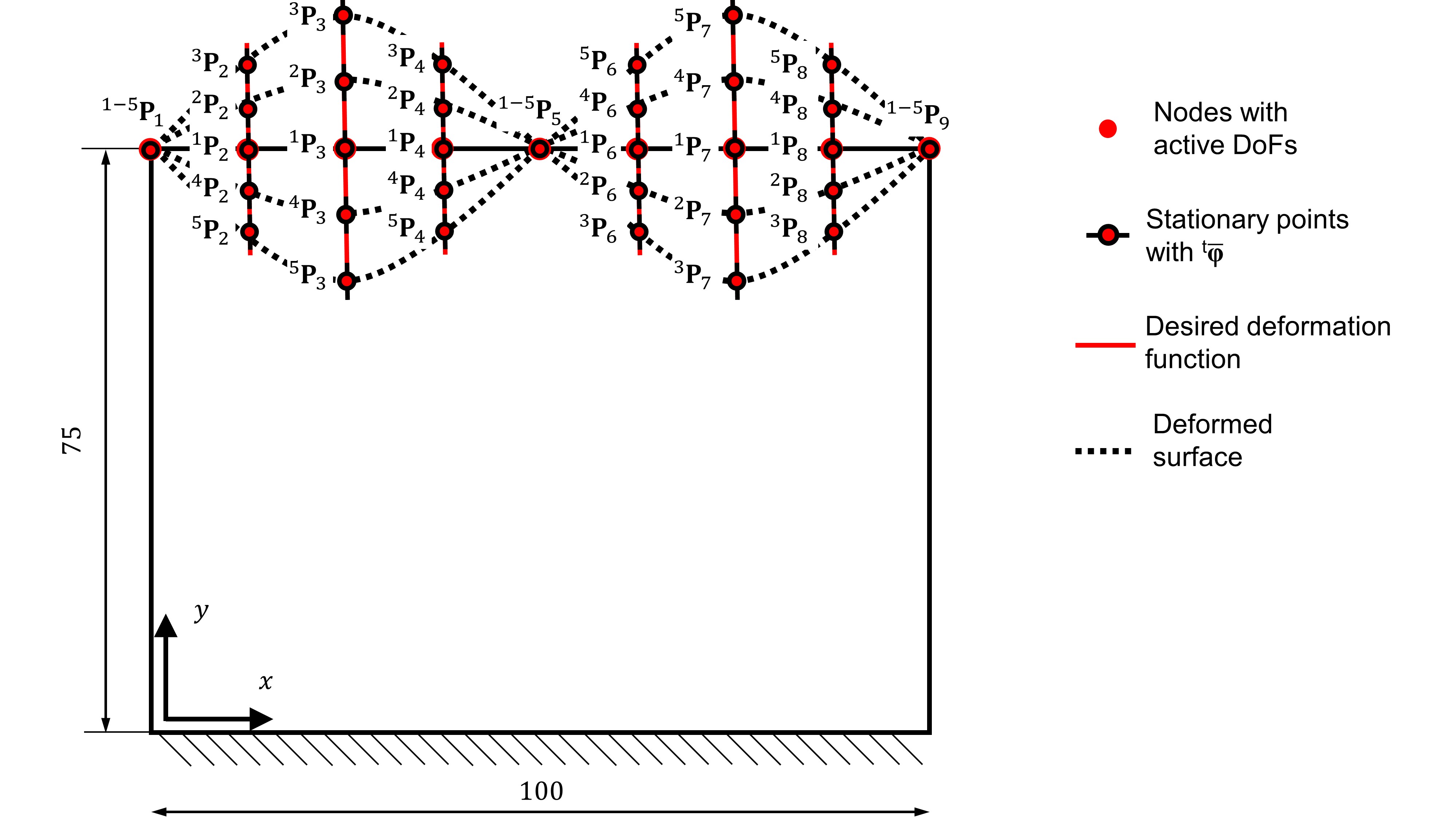}
    \caption{Shape-adaptive structure: design space with stationary points, desired deformation paths and the highlighted points with active DoFs \label{fig6}}
    \end{figure}

The design space shown in Figure~\ref{fig6} is used for the shape-adaptive structure. The design space is filled with 160 elements in the x-direction and 120 elements in the y-direction. The structure has 10 active DoFs. Nine DoFs represent a vertical sinusoidal displacement in the y-direction of the points highlighted in Figure~\ref{fig6}. In addition, the x-displacement of the point $^t\mathbf{P}_5$ is set to zero for all five stationary points. Except for point~5, the desired deformation paths are thus only defined in the y-direction. The $y$ displacements (referred to as $u_{ay}$) and the positions of the points in the $x$-direction are used for the plot of the displacements. The stationary points $^{2}\mathbf{u}_a$ and $^{4}\mathbf{u}_a$ describe a sine function with an amplitude of 1.25~mm and the stationary points $^{3}\mathbf{u}_a$ and $^{5}\mathbf{u}_a$ describe a sine function with an amplitude of 2.5~mm. 

    \unskip

The parameter $^t\upmu$ in \ref{Restriction_primary_stiffness2} is set to $^t\upmu=2.5$ for the linear and for the path-generating CM for each stationary point. For the path-generating shape-adaptive structure, the optimization was performed with variant~2 and $l$ is set to nine in \eqref{Restriction_mode_swap}. The weighting for the volume constraint in \ref{Max_secondary_stiffness_min_volume} is $\upzeta=0.7$ for both the linear and the path-generating CM.
\section{Results and Discussion}

\subsection{Performance parameters}

A CM performs well if the optimization goals (high accuracy and precision) are met as well as possible. Various characteristic values were introduced in the optimization of CM with linear FEM to investigate the performance of the designed CM. These can also be used in a modified form for path-generating CM. The ratio between the first two eigenvalues of the condensed stiffness matrix is referred to as selectivity $S$ and is a characteristic value for the precision of the CM. For high precision, this must be as high as possible \cite{Hasse2009}. The cosine similarity $\updelta$ was introduced as a characteristic value for accuracy in \cite{Kirmse2021}. This can assume values between zero and one. It is 1 if the two vectors are identical. The performance of the path-generating CM can be estimated if these two characteristic values are calculated for the tangent condensed stiffness matrices $^{t}\mathbf{\bar{K}}$ in the stationary points used for the optimization:

\begin{linenomath}
    \begin{equation}
        \label{selectivity}
        ^tS={^t\uplambda_{2}}/{^t\uplambda_{1}}   
    \end{equation}
\end{linenomath}

\begin{linenomath}
    \begin{equation}
        \label{cosine_similarity}
        ^t\updelta=\left| {^t\bar{\bm{\upchi}}_{1}^\mathrm{T}} {^t\bar{\bm{\upvarphi}}} \right|
    \end{equation}
\end{linenomath}

However, these values calculated for the stationary points only provide a good estimate if it is assumed that the stationary points are element of the natural kinematics.

For a more accurate but more complex evaluation, the performance of the CM can be determined in a modified form. If the CM is subjected to the respective first eigenmode in $^{t}\mathbf{u}_a$ step by step and new tangent stiffness matrices are then calculated, the natural kinematics as minimum energy path can be determined with sufficient accuracy if these steps are small enough:

\begin{linenomath}
    \begin{equation}
        \label{Kinematik}
        ^{t}\mathbf{\bar{K}}=\mathbf{\bar{K}}(^{t}\mathbf{u}_a),\;t=2...n_\mathrm{s}
    \end{equation}
\end{linenomath}

with

\begin{linenomath}
    \begin{equation}
        \begin{split}
        \label{Kinematik2}
        ^{1}\mathbf{u}_a&=\mathbf{0},\;t=1\\
        ^{t}\mathbf{u}_a&={^{t-1}\mathbf{u}_a}+\upbeta~{^t\bar{\bm{\upchi}}_{1}^\mathrm{T}},\;t=2...n_\mathrm{s}
        \end{split}
    \end{equation}
\end{linenomath}

The variable $n_\mathrm{s}$ is the number of calculated points for the natural kinematics. The natural kinematics are then compared with the selected stationary points for the optimization in order to estimate the accuracy. It is very high if the stationary points are element of the natural kinematics. The selectivity of the natural kinematics can also be calculated in the points of the natural kinematics in order to determine the precision. The precision can also be illustrated by loading the CM with selected load cases at the active DoFs and plotting the displacements in the active DoFs.

\subsection{Inverter}

\begin{figure}[H]
    \centering
    \includegraphics[width=0.8\textwidth]{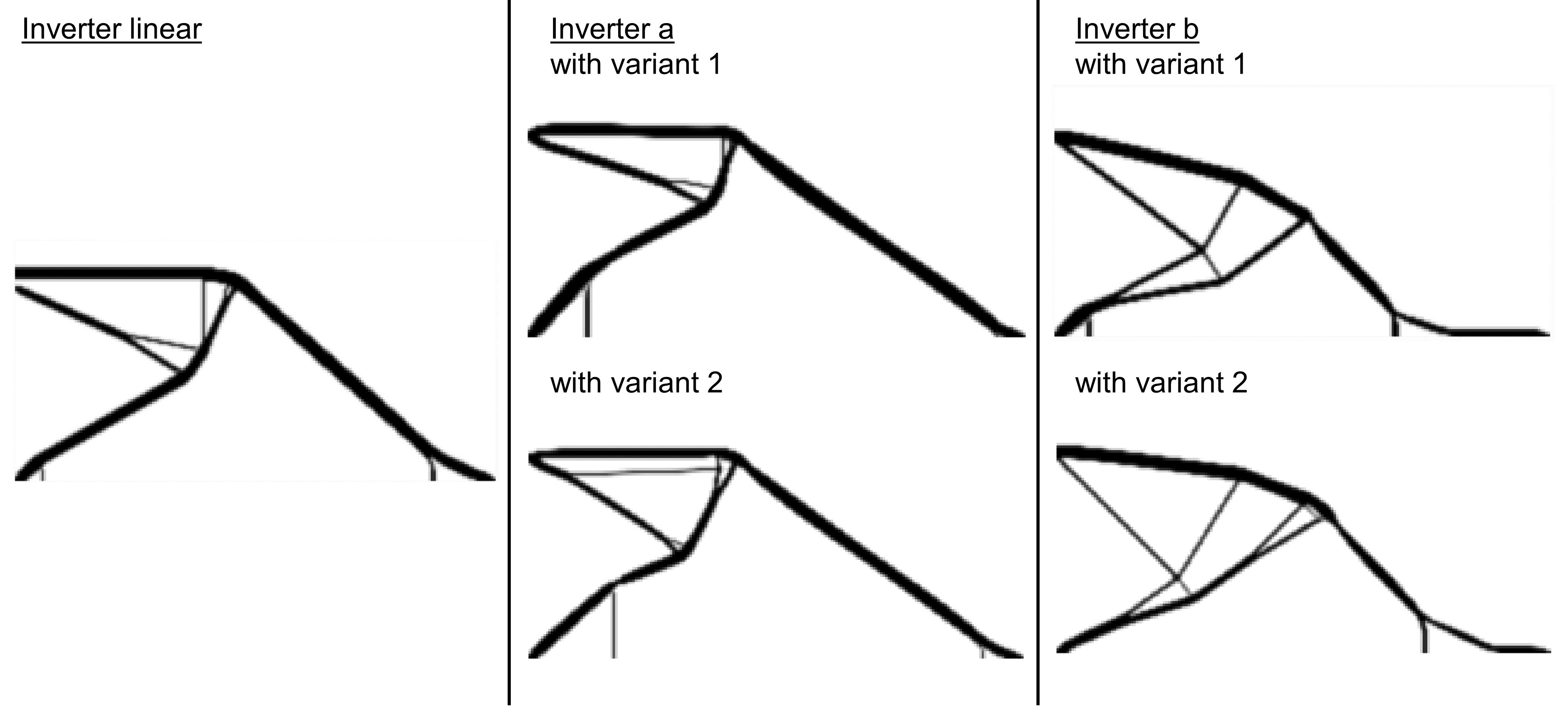}
    \caption{Inverter mechanism: topologies for various specified stationary points and optimization variants}
    \label{fig7}
\end{figure}

\begin{figure}[H]
    \centering
    \includegraphics[width=0.8\textwidth]{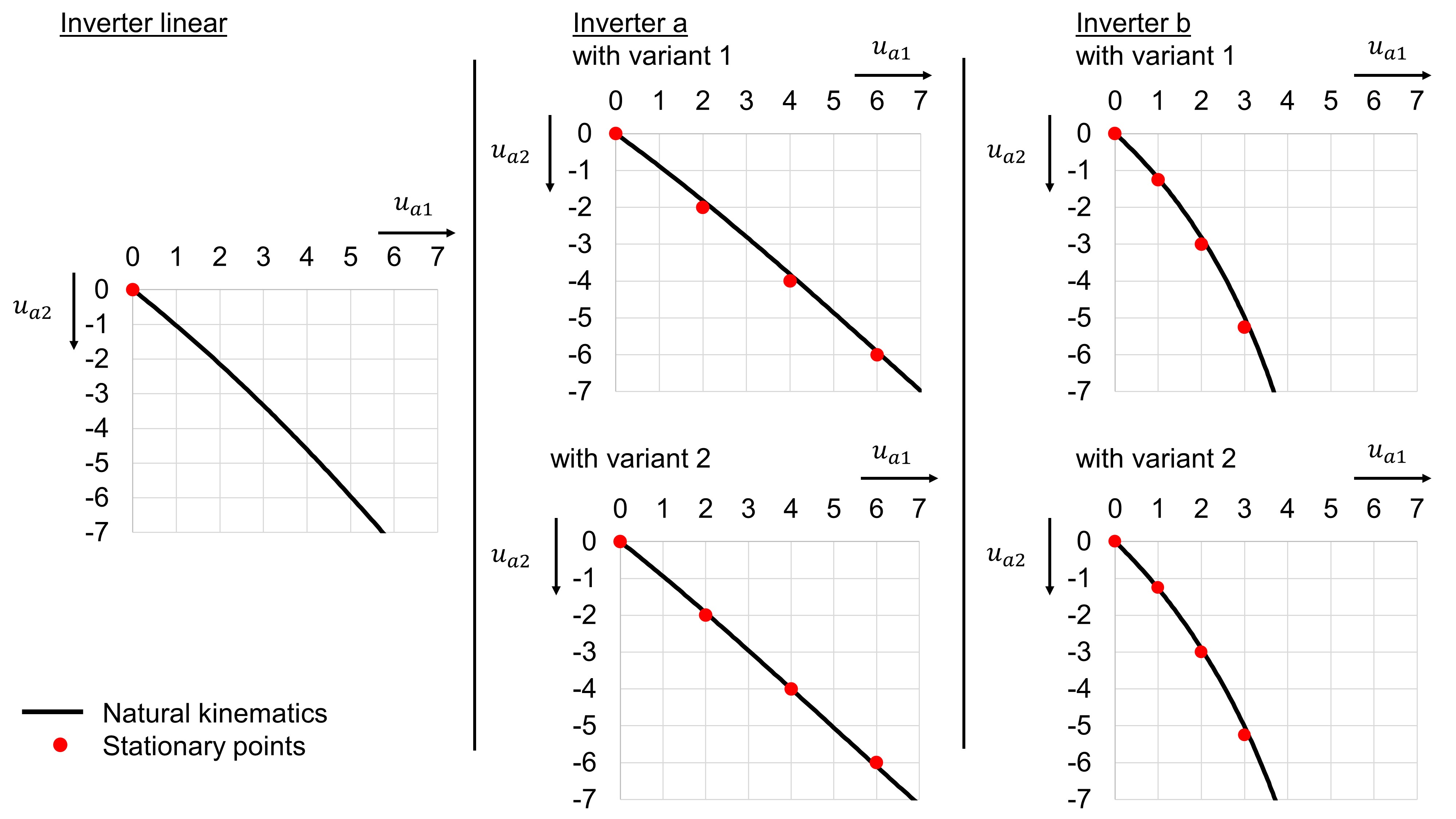}
    \caption{Inverter mechanism: natural kinematics for various specified stationary points and optimization variants}
    \label{fig8}
\end{figure}

The resulting topologies for the inverters are shown in Figure~\ref{fig7}. It can be seen that the inverters optimized with variant~2 have slightly thinner ribs than those optimized with variant~1. The reason for this is that the limit value for the primary stiffnesses can be better utilized with variant~1, as the primary stiffnesses can be freely distributed to the stationary points due to the summation function in equation \ref{Restriction_primary_stiffness1}. All path-generating inverters differ from the linear inverter, most clearly with inverter~b. To estimate the performance of these inverters, the selectivities and cosine similarities in the stationary points are summarized in Table~\ref{tab2}. It can be seen that the linear inverter has the highest selectivity at stationary point~1. However, the selectivities for inverter~a are also high in all other stationary points. The inverters~b have lower selectivities, particularly in the first stationary points. However, the selectivities increase sharply in the other stationary points. The cosine similarities are very high for all optimization variants and in all stationary points. The high selectivities and cosine similarities close to one indicate a good performance of the inverters. For a more precise statement, the natural kinematics are determined for all inverters for $\upbeta=0.1$. This is shown in Figure~\ref{fig8}. For all path-generating inverters, it can be seen that the natural kinematics approximate the specified stationary points well. Thus, good accuracy is achieved for all optimized inverters. The natural kinematics of the linear inverter is also shown. It can be seen that this differs from the two specified desired deformation functions for the path-generating inverters. It can be shown that different desired deformation functions can be specified for the optimization of path-generating inverters and that the natural kinematics are well approximated by the presented optimization method. 

\begin{figure}[H]
    \centering
    \includegraphics[width=0.8\textwidth]{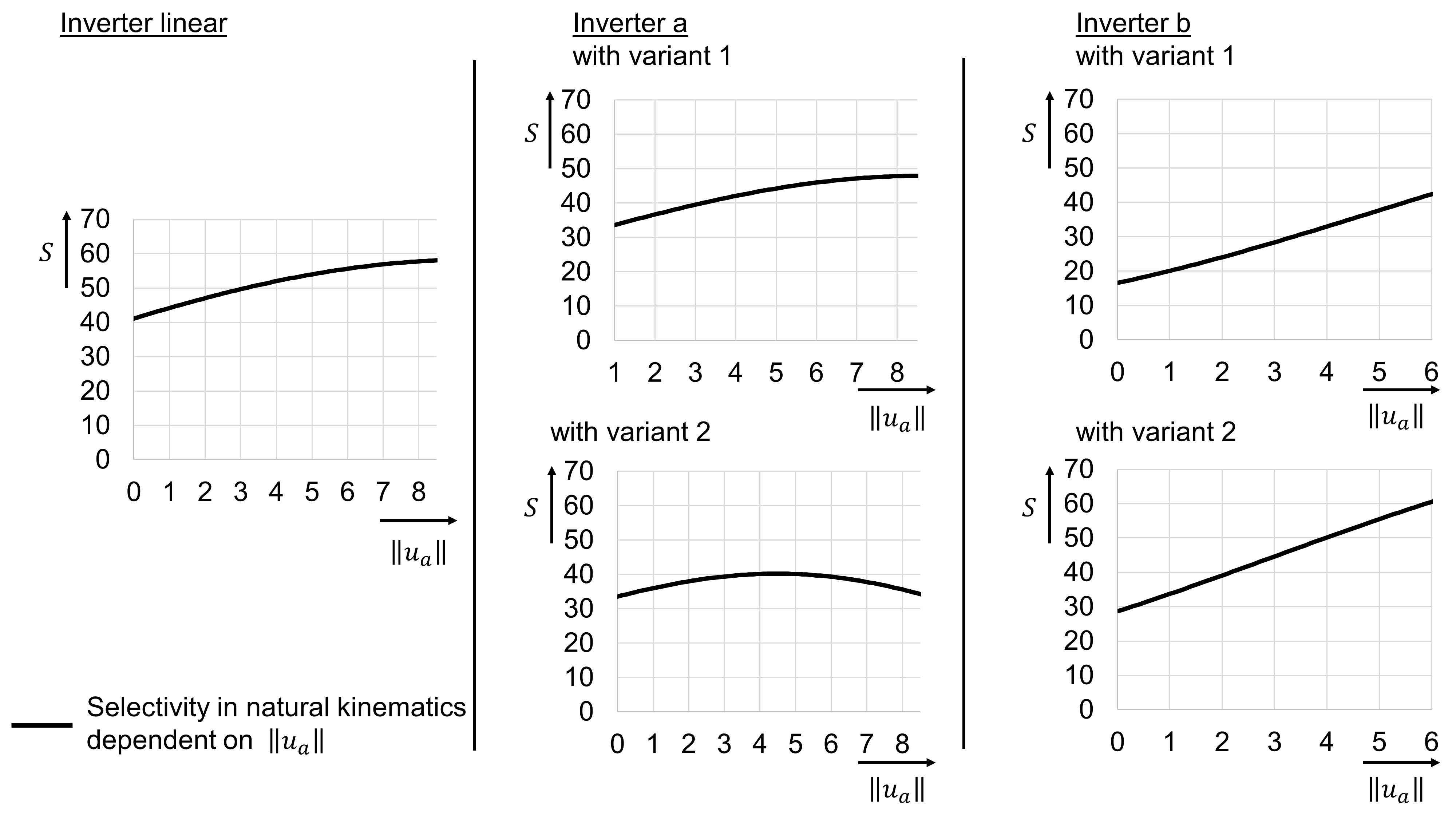}
    \caption{Inverter mechanism: selectivities of natural kinematics for various specified stationary points and optimization variants}
    \label{fig9}
\end{figure}

\begin{figure}[H]
    \centering
    \includegraphics[width=0.8\textwidth]{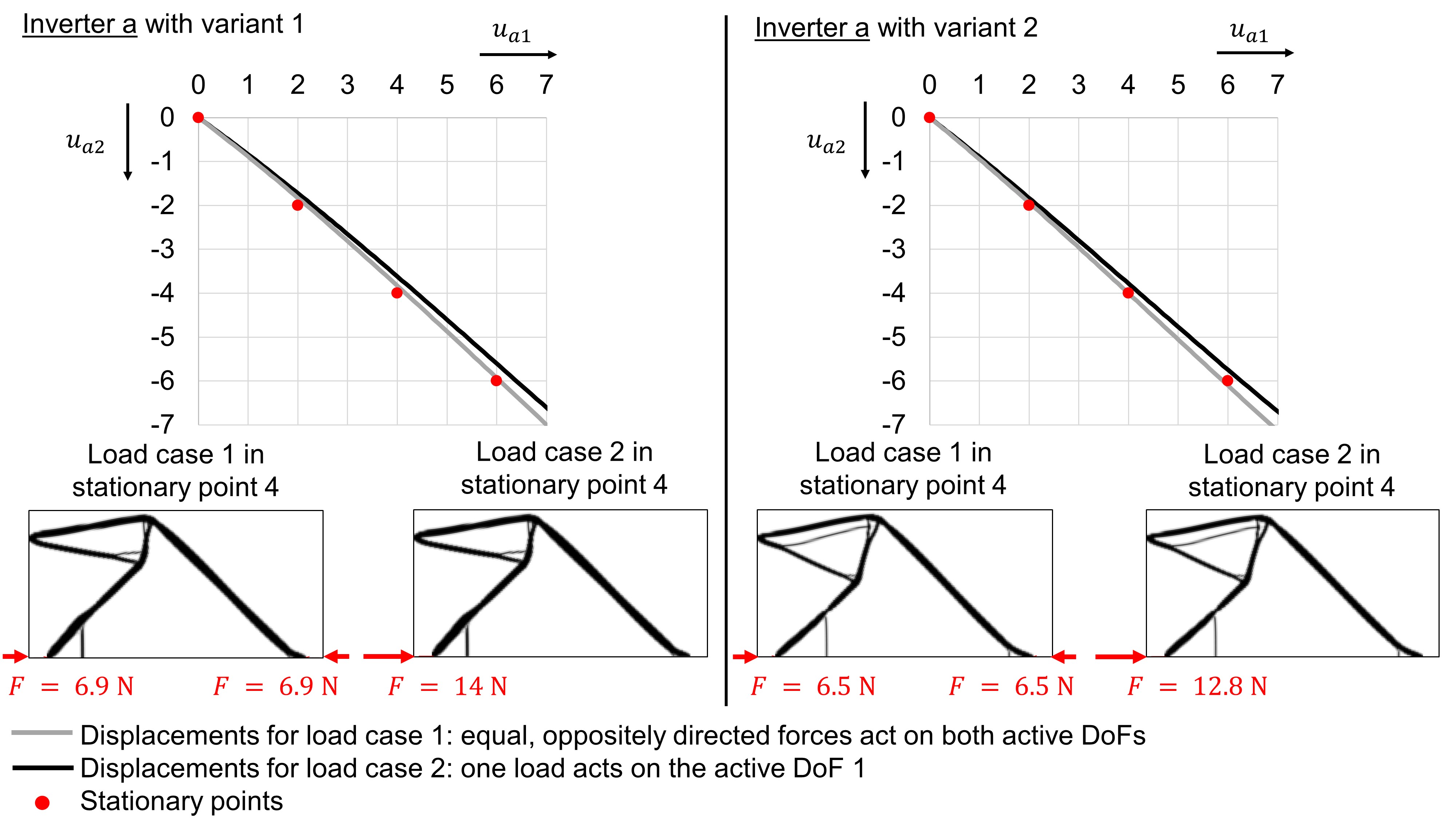}
    \caption{Inverter mechanism: displacement of the active DoFs when loading inverter~a with different load cases}
    \label{fig10}
\end{figure}

\begin{figure}[H]
    \centering
    \includegraphics[width=0.8\textwidth]{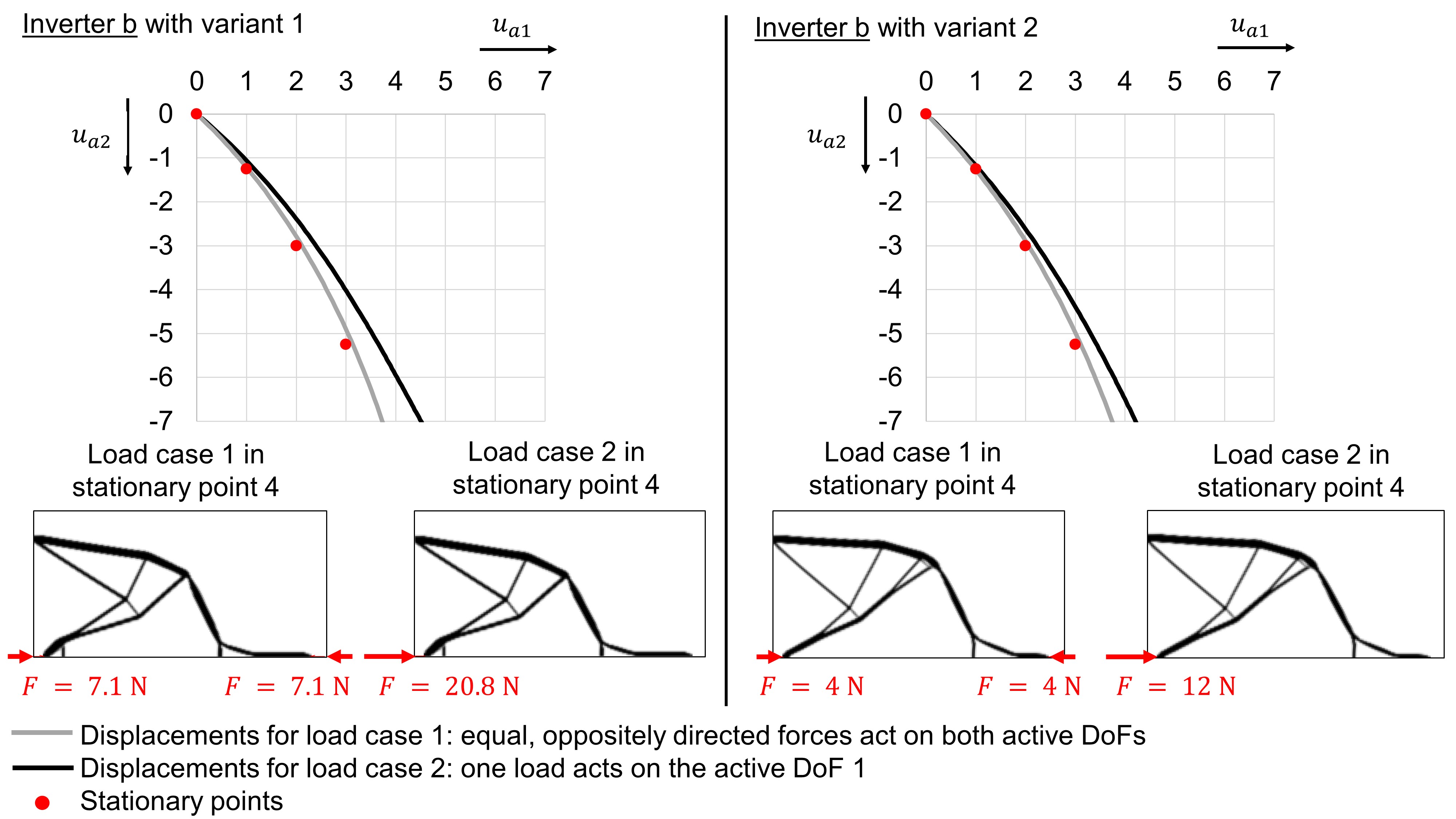}
    \caption{Inverter mechanism: displacement of the active DoFs when loading inverter~b with different load cases}
    \label{fig11}
\end{figure}

In order to be able to evaluate the precision, the selectivity of the natural kinematics is plotted in Figure~\ref{fig9} for the inverters as a function of the displacement $^t\mathbf{u}_{a1}$. In most cases, this increases with increasing displacement. For inverter~a in particular, the selectivities of the natural kinematics are very similar to the selectivities of the stationary points in Table~\ref{tab2}. For larger deformations (stationary point 3 and stationary point 4), however, the differences for inverter~b are larger. However, this does not affect the convergence of the optimization. The lower selectivities for inverter~b in the first stationary points suggest a more load-dependent kinematics than in the other examples. This is confirmed in Figure~\ref{fig10} and Figure~\ref{fig11}. It can be seen that inverter~a in \ref{fig10} has very good precision for various load cases. The ratio of the displacement between $^t\mathbf{u}_{a1}$ and $^t\mathbf{u}_{a2}$ remains almost the same regardless of the load. For inverter~b in \ref{fig11} the behavior is somewhat worse. Larger deviations occur with variant~1 in particular.

It can be summarized that the estimation of the characteristic values in the stationary points basically provides a good indication of the performance of the CM. For more precise information the natural kinematics should be determined. It has been shown that it is possible to design inverters that follow very different deformation paths. Good accuracy can also be guaranteed. No recommendation can be given as to which variant should be used for optimization. Both variants generally deliver good results and are therefore well suited.

\unskip

\begin{table}[H] 
    \caption{Selectivity and cosine similarity for the inverter mechanisms\label{tab2}}
    \newcolumntype{C}{>{\centering\arraybackslash}X}
    \begin{tabularx}{\textwidth}{lC|CCCCCCC}
    \toprule
    ~ & Stationary point & Linear & Inverter a & Inverter a & Inverter b & Inverter b\\
    ~ & ~ & ~ & Variant 1 & Variant 2 & Variant 1 & Variant 2 \\
    \midrule
    {\multirow{5}{*}{$^tS$}} 
    & 1 & 41.1365 & 30.5009 & 33.5475 & 16.6031 & 28.6846 \\
    & 2 & - & 42.5866 & 40.0940 & 23.2096 & 35.7115\\
    & 3 & - & 47.7834 & 39.8668 & 36.8087 & 54.1719 \\
    & 4 & - & 48.1845 & 34.7708 & 52.1825 & 77.5740 \\
    \midrule
    {\multirow{5}{*}{$^t\updelta$}} 
    & 1 & 1.0000 & 0.9970 & 0.9990 & 0.9995 & 0.9979 \\
    & 2 & - & 0.9999 & 1.0000 & 0.9995 & 0.9998 \\
    & 3 & - & 0.9998 & 0.9998 & 0.9997 & 0.9996 \\
    & 4 & - & 0.9994 & 0.9998 & 0.9998 & 1.0000 \\
    \bottomrule
    \end{tabularx}
    \end{table}
    \unskip

\unskip

\unskip

\unskip

\unskip

\subsection{Pivot joint}

\begin{figure}[H]
    \centering
    \includegraphics[width=0.8\textwidth]{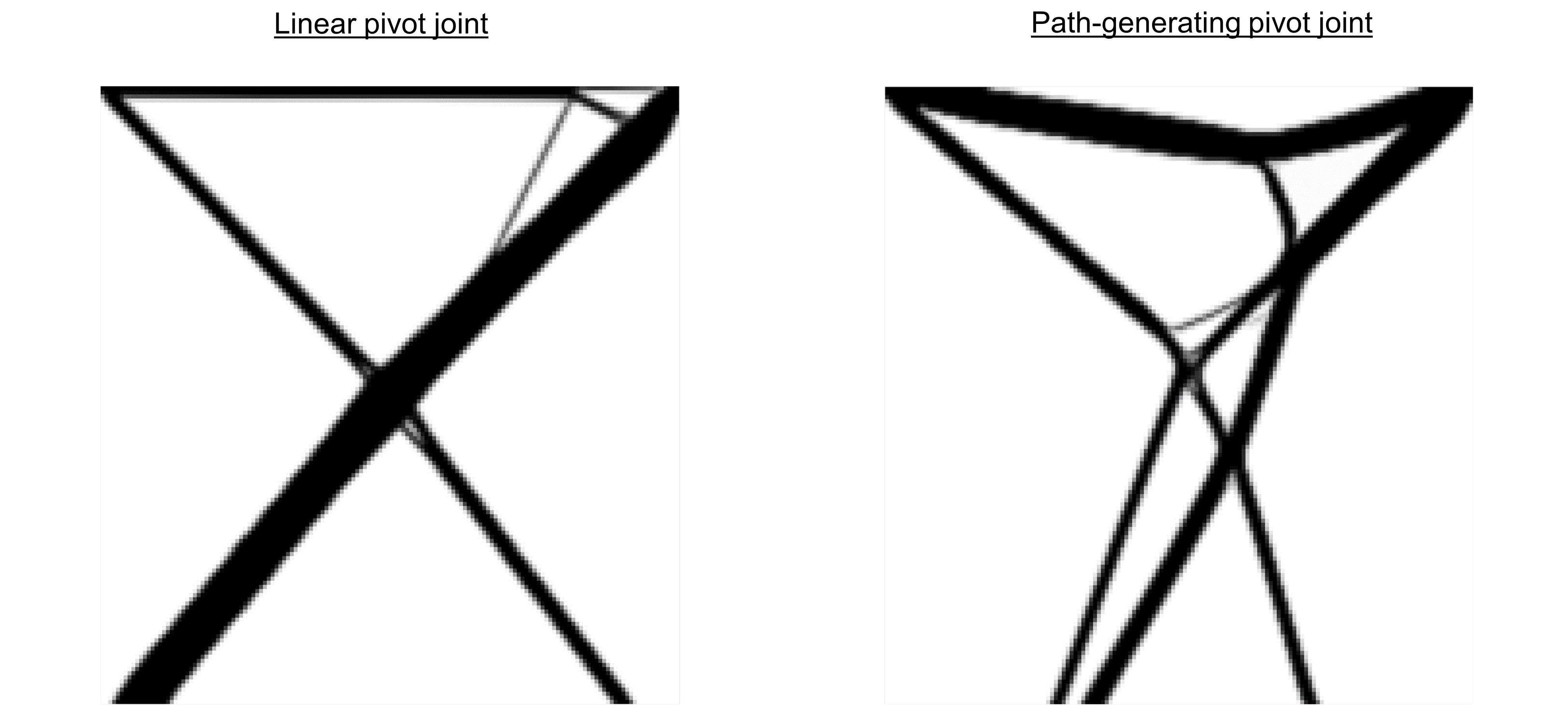}
    \caption{Pivot joint: topologies for various specified stationary points}
    \label{fig12}
\end{figure}

\begin{figure}[H]
    \centering
    \includegraphics[width=0.8\textwidth]{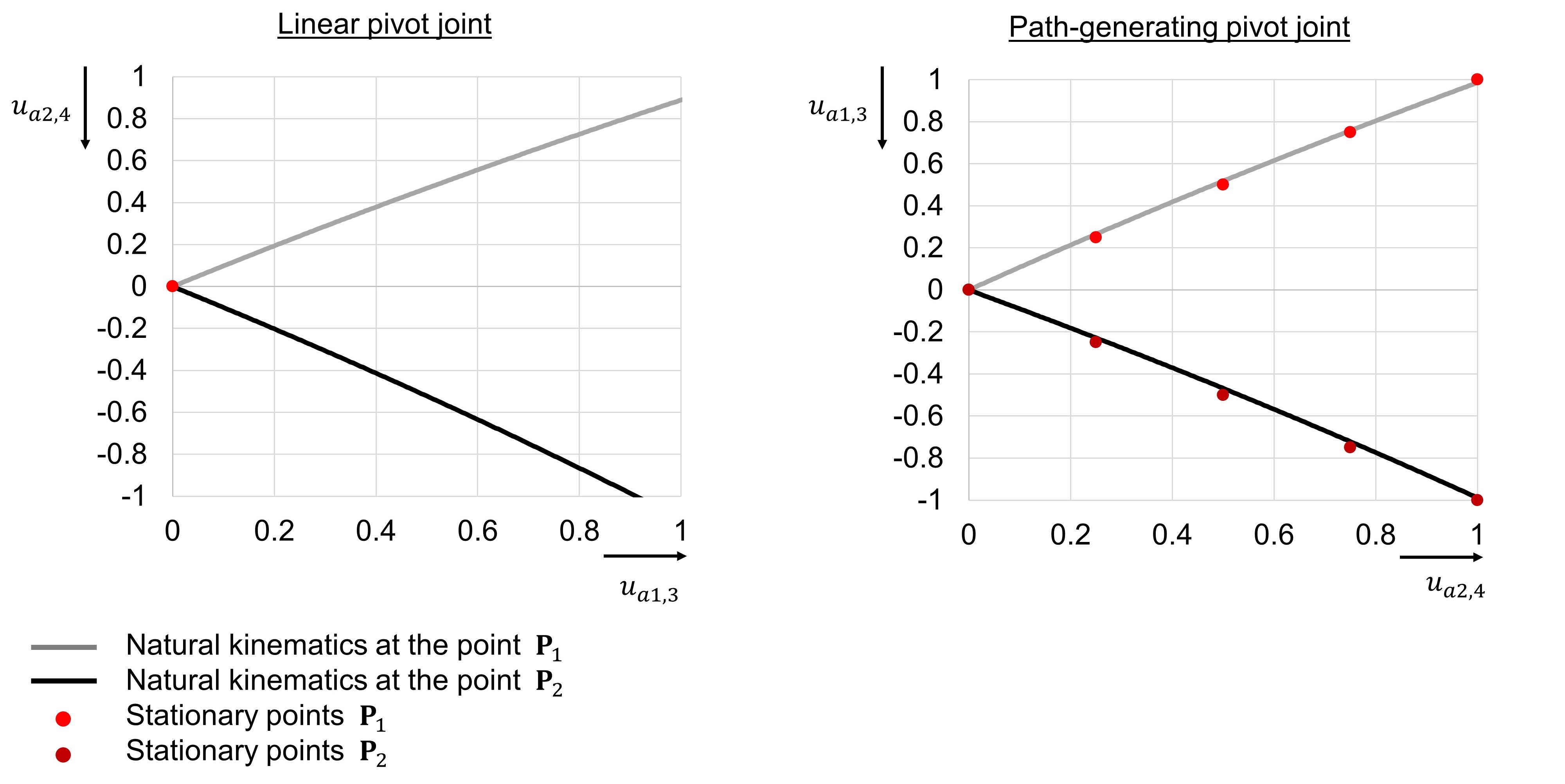}
    \caption{Pivot joint: natural kinematics for various specified stationary points}
    \label{fig13}
\end{figure}

\begin{figure}[H]
    \centering
    \includegraphics[width=0.8\textwidth]{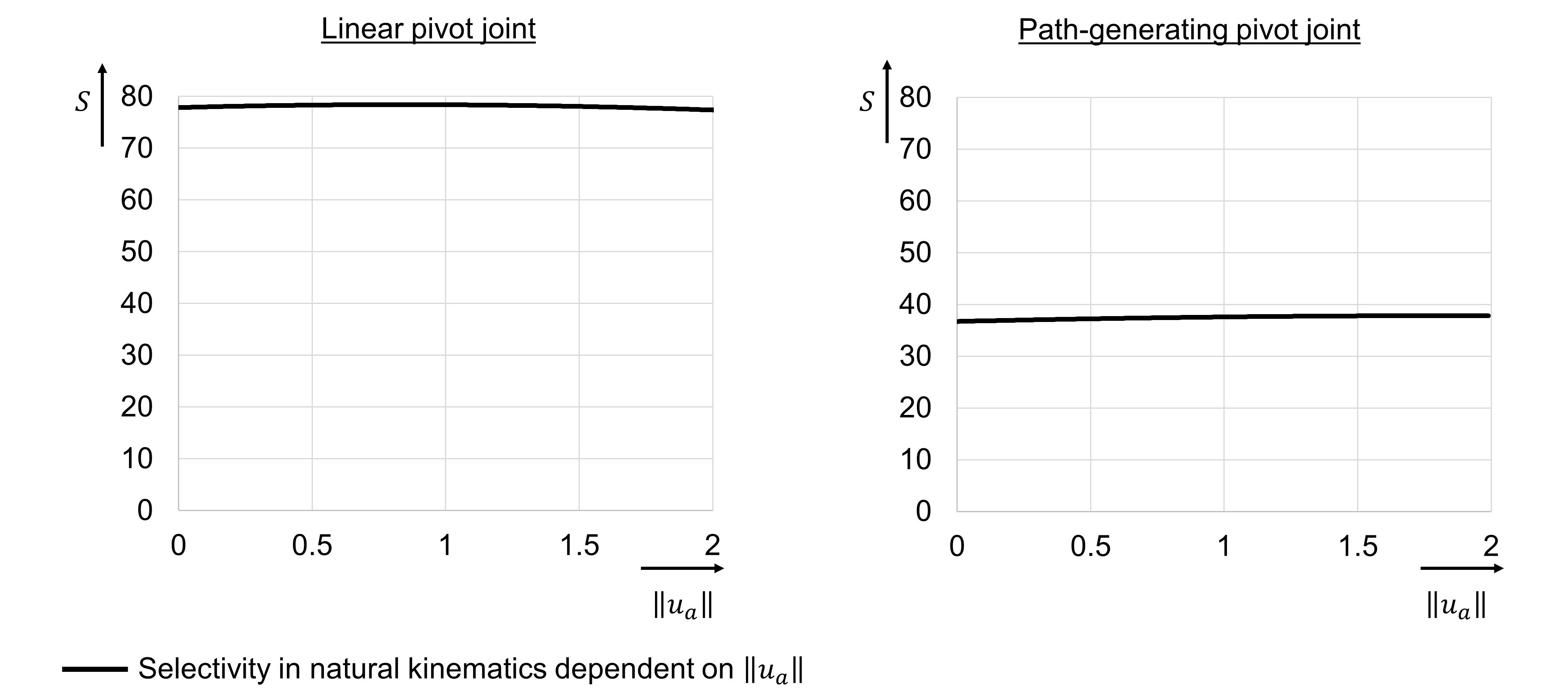}
    \caption{Pivot joint: selectivities of the natural kinematics for various specified stationary points}
    \label{fig14}
\end{figure}

\begin{figure}[H]
    \centering
    \includegraphics[width=0.8\textwidth]{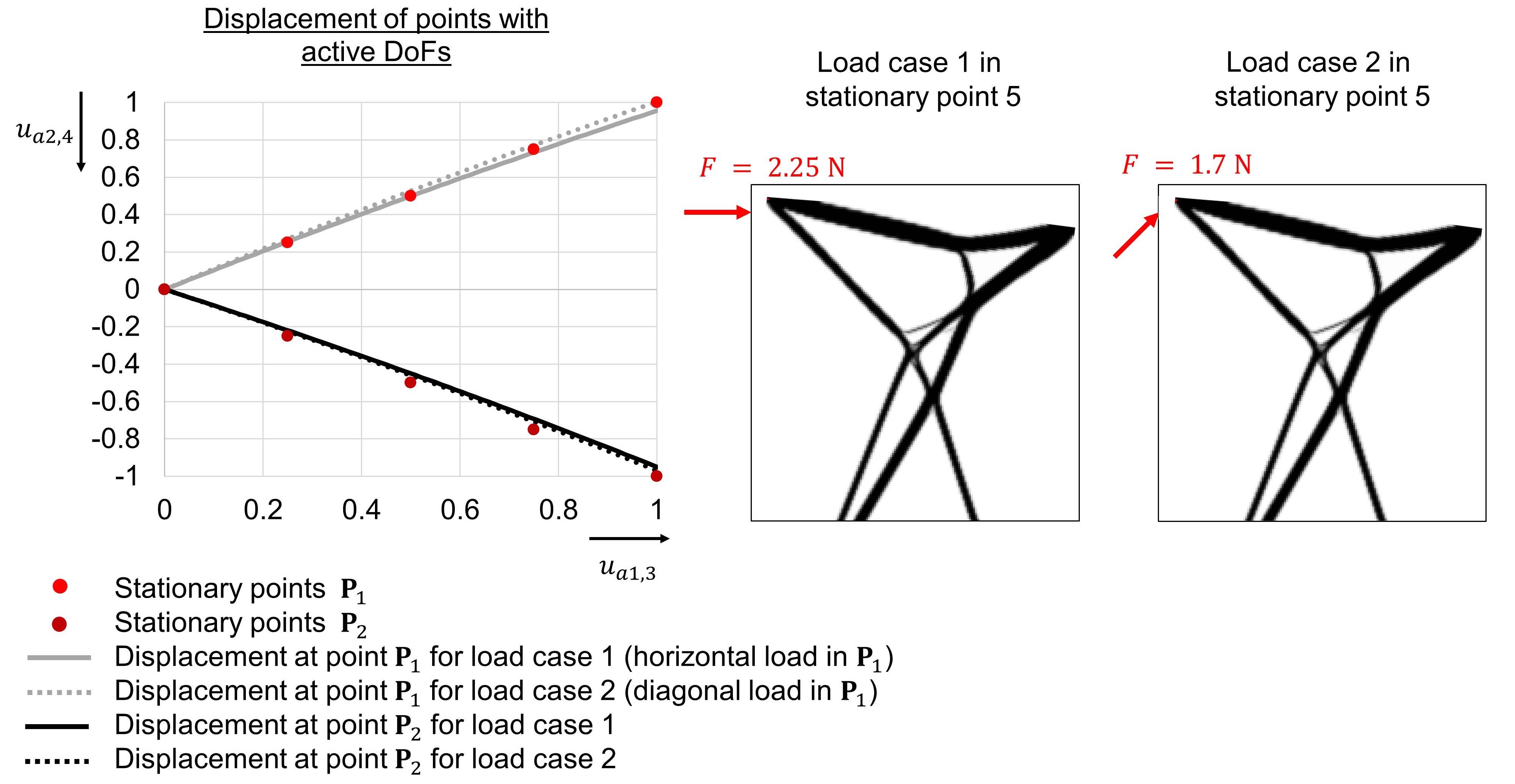}
    \caption{Displacement of the active DoFs when loading the path-generating pivot joints with different load cases}
    \label{fig15}
\end{figure}

The optimized pivot joints are shown in Figure~\ref{fig12}. It can be seen that the optimized structures differ greatly. While the linear pivot joint is similar to the familiar cartwheel joint, the structure for the path-generating CM is more complex. The performance parameters in the stationary points are summarized in Table~\ref{tab3}. It can be seen that high cosine similarities and selectivities are achieved for both CMs in all stationary points. The natural kinematics of the two structures for $\upbeta=0.01$ is shown in Figure~\ref{fig13}. It can be seen that the natural kinematics are well approximated to the specified stationary points by the synthesis as a path-generating CM. The selectivity of the natural kinematics shown in Figure~\ref{fig14} is high for both structures with increasing deformation. The precision achieved with path-generating CM is so high that there are hardly any load-dependent deviations in the deformation of the active DoFs, as shown in Figure~\ref{fig15} for two selected load cases. 

\unskip

\begin{table}[H] 
    \caption{Selectivity and cosine similarity for the pivot joints\label{tab3}}
    \newcolumntype{C}{>{\centering\arraybackslash}X}
    \begin{tabularx}{\textwidth}{C|CCCCC}
    \toprule
    {\multirow{2}{*}{\makecell[tc]{Stationary \\ point}}}
    & $^tS$ & $^tS$  & $^t\updelta$ & $^t\updelta$ \\
    & Linear & Path-generating & Linear & Path-generating \\
    \midrule
    1 & 77.8589 & 36.7245 & 1.0000 & 0.9986 \\
    2 & - & 59.3637 & - & 0.9997 \\
    3 & - & 71.8011 & - & 1.0000 \\
    4 & - & 53.1275 & - & 0.9996 \\
    5 & - & 32.6779 & - & 0.9984 \\
    \bottomrule
    \end{tabularx}
    \end{table}
    \unskip

\unskip

\unskip

\unskip

\subsection{Shape-adaptive structure}

\begin{figure}[H]
    \centering
    \includegraphics[width=0.8\textwidth]{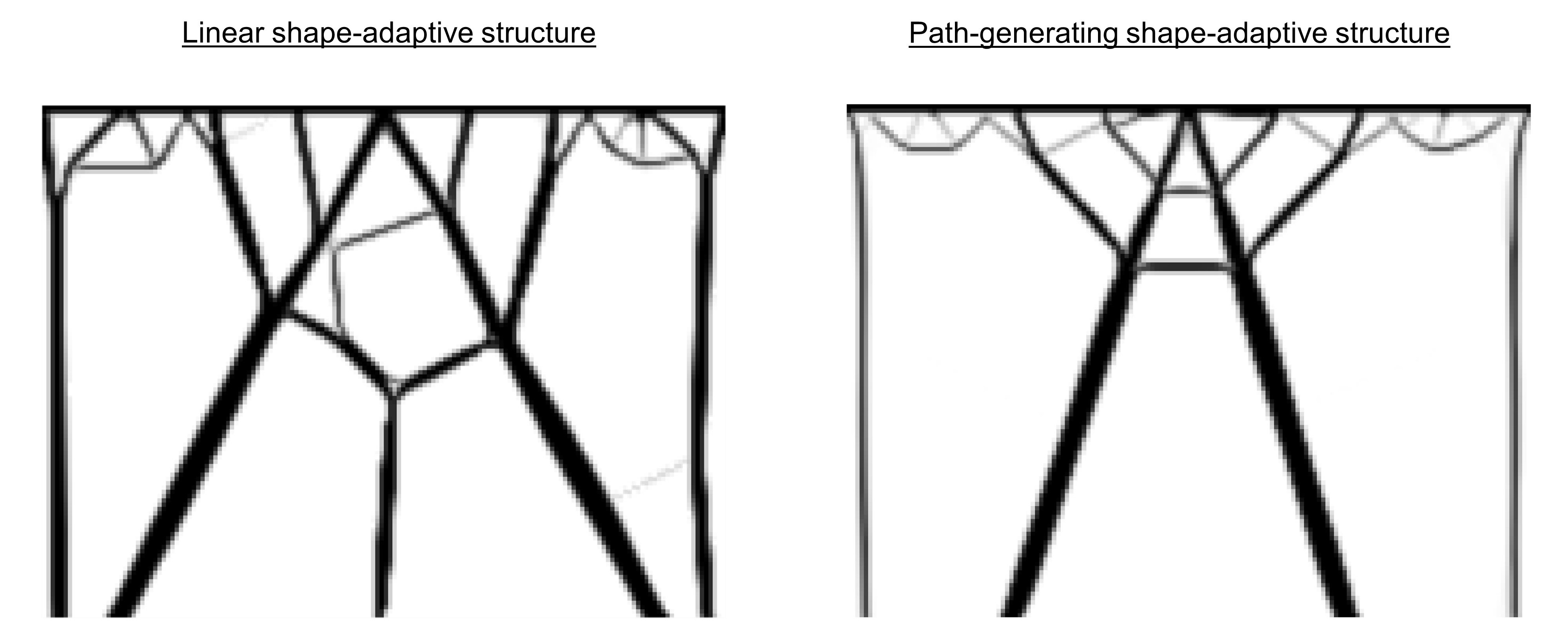}
    \caption{Shape-adaptive structure: topologies for various specified stationary points}
    \label{fig19}
\end{figure}

\begin{figure}[H]
        \centering
        \includegraphics[width=0.8\textwidth]{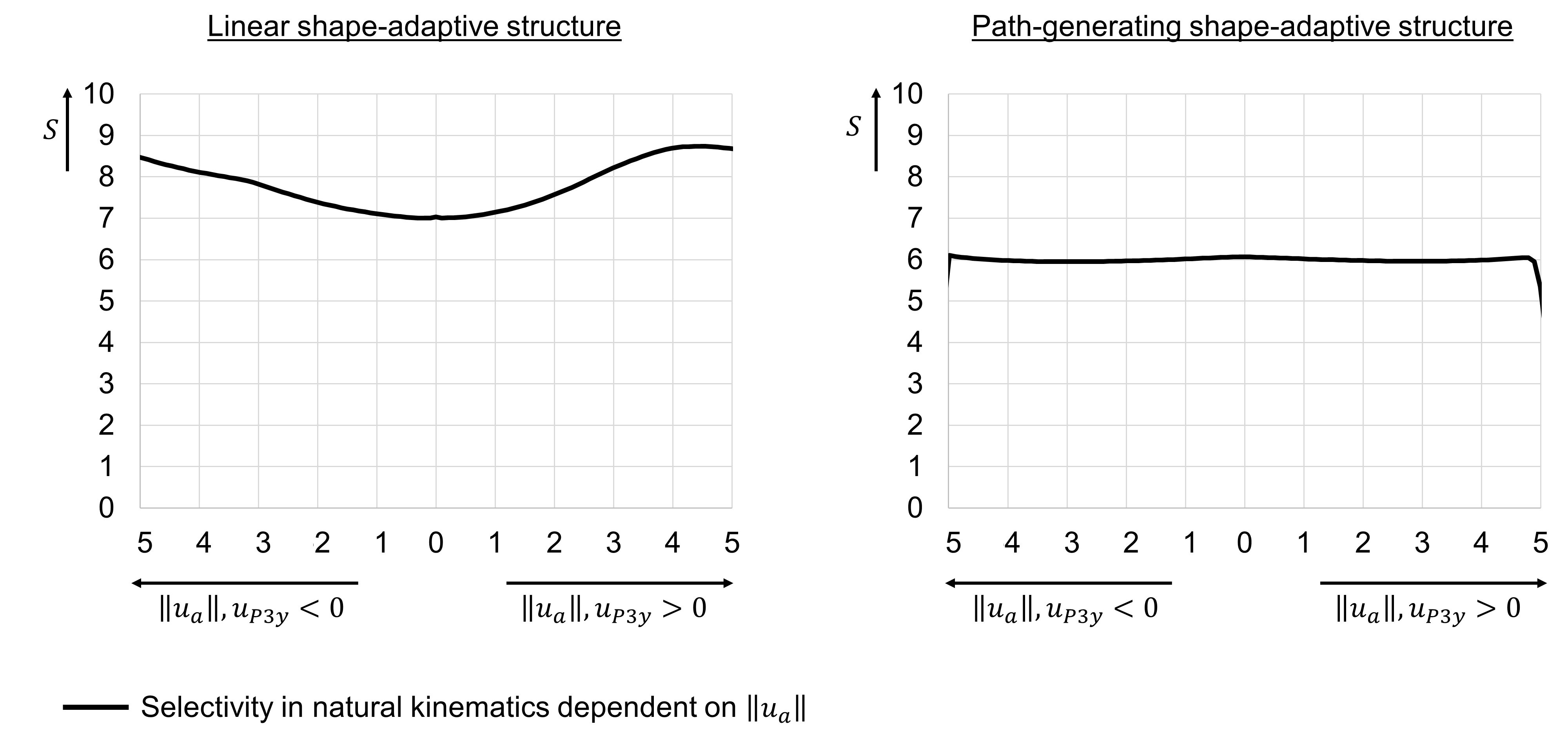}
        \caption{Shape-adaptive structure: selectivities of the natural kinematics for various specified stationary points}
        \label{fig20}
    \end{figure}

\begin{figure}[H]
        \centering
        \includegraphics[width=0.8\textwidth]{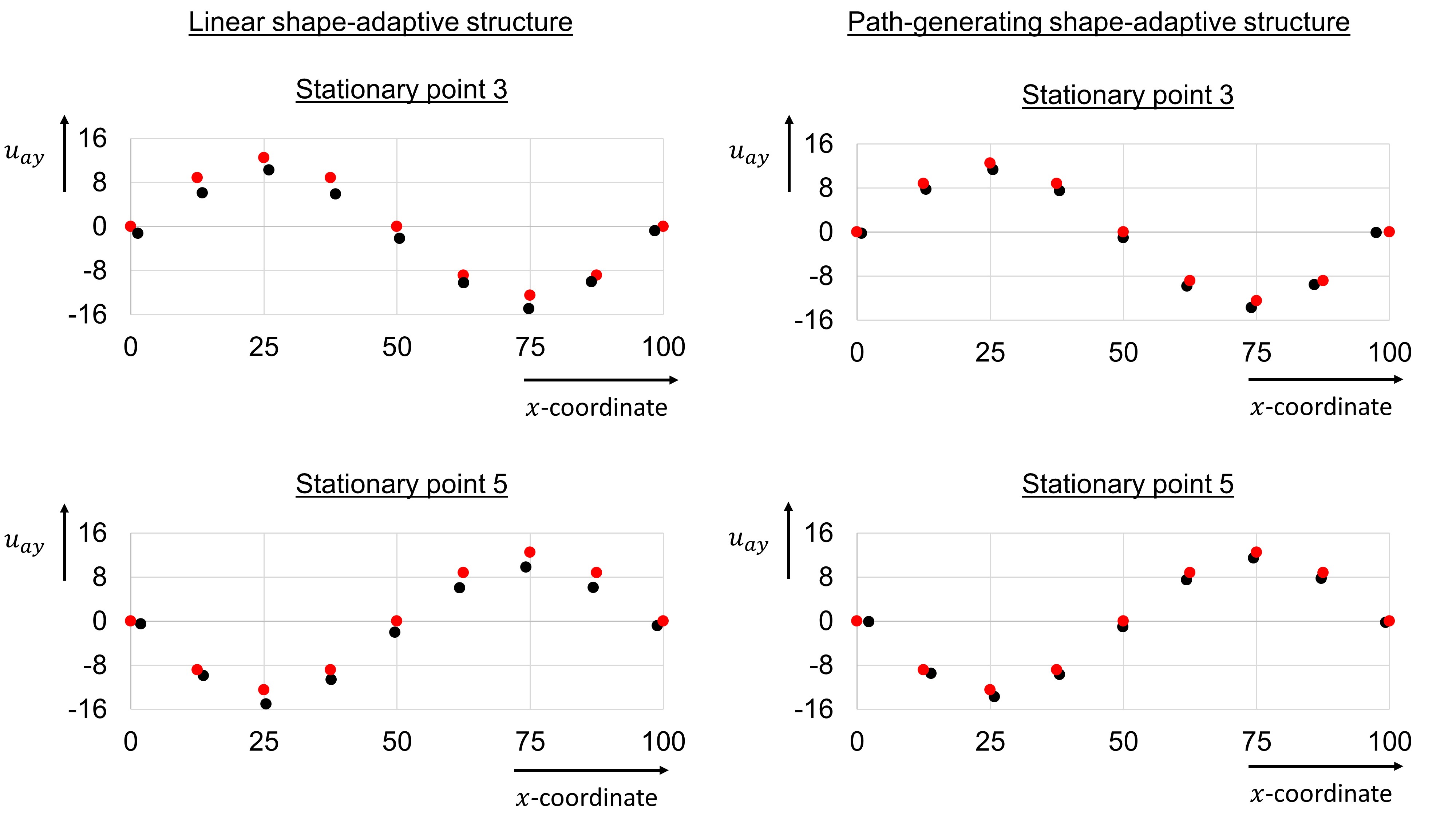}
        \caption{Shape-adaptive structure: natural kinematics for various specified stationary points, scaled by a factor of 5}
        \label{fig21}
    \end{figure}

\begin{figure}[H]
        \centering
        \includegraphics[width=0.8\textwidth]{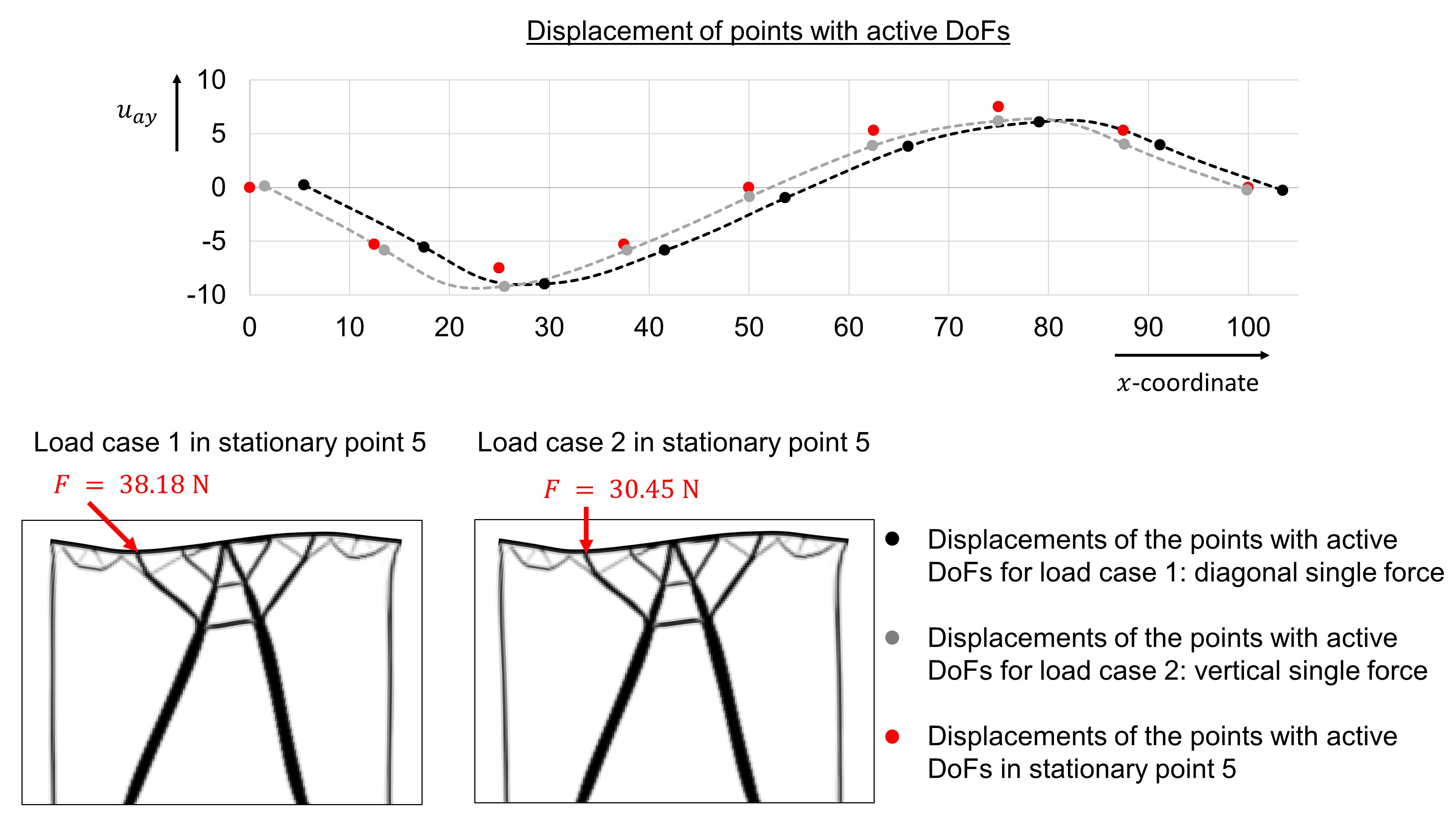}
        \caption{Displacement of the active DoFs when the path-generating shape-adaptive structure is loaded with different load cases, scaled by a factor of 3}
        \label{fig22}
    \end{figure}

The optimization results for the shape-adaptive structures are shown in Figure~\ref{fig19}. The optimized structures also differ greatly from each other in this design example. The performance parameters in the stationary points are listed in Table~\ref{tab4}. In the first stationary point, the selectivity is higher for the linear structure. However, the path-generating structure has a better cosine similarity in this stationary point. At the other stationary points, the path-generating structure exhibits high cosine similarities and high selectivities. It should be noted that the selectivities here are lower than in the previous design examples. However, the deformation is more complex. The selectivity in natural kinematics for $\upbeta=0.1$ is shown in Figure~\ref{fig20}. It is slightly lower for the path-generating structure than for the linear structure. It can be seen that the selectivity decreases at the edges of the domain of definition. This is due to the fact that the mechanism buckles. However, buckling does not occur in the domain of definition. The natural kinematics is shown in Figure~\ref{fig21} for the case that the magnitude of $\mathbf{u}_a$ is equal to the magnitude of the shown stationary points. It can be seen that a better approximation can be achieved with the path-generating structure than with the linear structure. It can therefore be shown that a path-generating shape-adaptive structure with selective compliance can be designed using the presented optimization approach. The selectively compliant behavior is shown in Figure~\ref{fig22}. For this purpose, the path-generating shape-adaptive structure is loaded with different forces and the deformation is shown with a similar amount in the active DoFs as in stationary point 5. The load-dependent deviations are small.

\unskip

\begin{table}[H] 
    \caption{Selectivity and cosine similarity for the shape-adaptive structure \label{tab4}}
    \newcolumntype{C}{>{\centering\arraybackslash}X}
    \begin{tabularx}{\textwidth}{C|CC|CC}
    \toprule
    {\multirow{2}{*}{\makecell[tc]{Stationary \\ point}}}
    & $^tS$ & $^tS$  & $^t\updelta$ & $^t\updelta$ \\
    & Linear & Path-generating & Linear & Path-generating \\
    \midrule
    1 &   7.0333   &  6.0701 & 0.9982 & 0.9998 \\
    2 & - &  5.6804 & - & 0.9929 \\
    3 & - & 5.0186 & - & 0.9719 \\
    4 & - & 5.6797 & - & 0.9933 \\
    5 & - & 4.9710 & - & 0.9725 \\
    \bottomrule
    \end{tabularx}
    \end{table}
    \unskip

    \unskip
    
    \unskip
    
    \unskip
\section{Conclusions}

The research presented in this paper addresses the limitations of current synthesis approaches for compliant mechanisms (CMs) by extending the topology optimization process to account for geometric and material nonlinearities. Traditional methods often rely on assumptions of small distortions, which constrain their applicability to real-world scenarios where CMs typically experience large deformations. As a result, existing methods have been inadequate for ensuring reliable path-generating behavior under varied load conditions.

The authors' prior work on developing load case insensitive CMs laid the groundwork for this study. By extending that pseudo-kinematic approach to include nonlinearities, the research advances the synthesis of path-generating CMs capable of maintaining their designed deformation paths even under diverse loading conditions. This improvement is crucial for the practical deployment of CMs in complex applications.

The effectiveness of the proposed method has been validated through several design examples, demonstrating its robustness and versatility. Furthermore, this paper introduces a novel shape-adaptive path-generating CM, showcasing the potential for creating more sophisticated and adaptable mechanisms.

In summary, this research contributes significant advancements to the field of compliant mechanism design by overcoming previous limitations related to nonlinearity and load sensitivity. The extended approach offers a more reliable and practical solution for synthesizing CMs, paving the way for their broader application in engineering and technology. Future work will focus on enhancing selectivity and modeling efficiency, as well as designing CMs capable of handling larger deformations.

\section{Declaration of Competing Interest}

The authors declare that they have no known competing financial interests or personal relationships that could have appeared to influence the work reported in this paper.

\section{Funding}

This work was funded by the Deutsche Forschungsgemeinschaft (DFG, German Research Foundation), project number HA 7893/3-2.

\section{Acknowledgement}

The authors gratefully acknowledge the GWK support for funding this project by providing computing time through the Center for Information Services and HPC (ZIH) at TU Dresden.

\printbibliography

@book{Bathe2014,
  title = {Finite Element Procedures},
  author = {Bathe, Klaus-J{\"u}rgen},
  year = {2014},
  edition = {2nd ed},
  publisher = {Prentice-Hall},
  address = {Englewood Cliffs, N.J},
  isbn = {978-0-9790049-5-7},
  langid = {english}
}

@incollection{Bendsoe2003c,
  title = {Design Parametrization},
  booktitle = {Topology Optimization: Theory, Methods, and Applications},
  author = {Bends{\o}e, Martin P. and Sigmund, O.},
  year = {2003},
  pages = {4--8},
  publisher = {Springer},
  address = {Berlin ; New York},
  isbn = {978-3-540-42992-0},
  langid = {english},
  lccn = {TA658.8 .B463 2003}
}

@incollection{Bendsoe2003d,
  title = {Problem Setting},
  booktitle = {Topology Optimization: Theory, Methods, and Applications},
  author = {Bends{\o}e, Martin P. and Sigmund, O.},
  year = {2003},
  pages = {95--97},
  publisher = {Springer},
  address = {Berlin ; New York},
  isbn = {978-3-540-42992-0},
  langid = {english},
  lccn = {TA658.8 .B463 2003}
}

@article{Bourdin2001,
  title = {Filters in Topology Optimization},
  author = {Bourdin, Blaise},
  year = {2001},
  journal = {International Journal for Numerical Methods in Engineering},
  volume = {50},
  number = {9},
  pages = {2143--2158},
  issn = {0029-5981, 1097-0207},
  doi = {10.1002/nme.116},
  urldate = {2022-11-03},
  langid = {english}
}

@inproceedings{Bruns1998,
  title = {Topology Optimization of Geometrically Nonlinear Structures and Compliant Mechanisms},
  booktitle = {7th {{AIAA}}/{{USAF}}/{{NASA}}/{{ISSMO Symposium}} on {{Multidisciplinary Analysis}} and {{Optimization}}},
  author = {Bruns, Tyler and Tortorelli, Daniel},
  year = {1998},
  pages = {4950},
  publisher = {{American Institute of Aeronautics and Astronautics}},
  address = {St. Louis,MO,U.S.A.},
  doi = {10.2514/6.1998-4950},
  urldate = {2020-04-29},
  langid = {english}
}

@article{Bruns2001a,
  title = {Topology Optimization of Non-Linear Elastic Structures and Compliant Mechanisms},
  author = {Bruns, Tyler E. and Tortorelli, Daniel A.},
  year = {2001},
  journal = {Computer Methods in Applied Mechanics and Engineering},
  volume = {190},
  number = {26-27},
  pages = {3443--3459},
  issn = {00457825},
  doi = {10.1016/S0045-7825(00)00278-4},
  urldate = {2022-11-03},
  langid = {english}
}

@article{Bruns2003,
  title = {An Element Removal and Reintroduction Strategy for the Topology Optimization of Structures and Compliant Mechanisms},
  author = {Bruns, T. E. and Tortorelli, D. A.},
  year = {2003},
  journal = {International Journal for Numerical Methods in Engineering},
  volume = {57},
  number = {10},
  pages = {1413--1430},
  issn = {0029-5981, 1097-0207},
  doi = {10.1002/nme.783},
  urldate = {2020-04-29},
  langid = {english}
}

@article{Bruns2005,
  title = {A Reevaluation of the {{SIMP}} Method with Filtering and an Alternative Formulation for Solid--Void Topology Optimization},
  author = {Bruns, T.E.},
  year = {2005},
  journal = {Structural and Multidisciplinary Optimization},
  volume = {30},
  number = {6},
  pages = {428--436},
  issn = {1615-147X, 1615-1488},
  doi = {10.1007/s00158-005-0537-x},
  urldate = {2022-02-07},
  langid = {english}
}

@article{Campanile2022,
  title = {Accuracy and Precision: {{A}} New View on Kinematic Assessment of Solid-State Hinges and Compliant Mechanisms},
  shorttitle = {Accuracy and Precision},
  author = {Campanile, Lucio Flavio and Kirmse, Stephanie and Hasse, Alexander},
  year = {2022},
  journal = {Journal of Intelligent Material Systems and Structures},
  volume = {33},
  number = {13},
  pages = {1743--1748},
  issn = {1045-389X, 1530-8138},
  doi = {10.1177/1045389X211064323},
  urldate = {2023-03-02},
  langid = {english}
}

@article{Campello2003,
  title = {A Triangular Finite Shell Element Based on a Fully Nonlinear Shell Formulation},
  author = {Campello, E. M. B. and Pimenta, P. M. and Wriggers, P.},
  year = {2003},
  journal = {Computational Mechanics},
  volume = {31},
  number = {6},
  pages = {505--518},
  issn = {0178-7675, 1432-0924},
  doi = {10.1007/s00466-003-0458-8},
  urldate = {2023-02-28},
  langid = {english}
}

@article{Capasso2020,
  title = {Stress-Based Topology Optimization of Compliant Mechanisms Using Nonlinear Mechanics},
  author = {Capasso, Gabriele and Morlier, Joseph and Charlotte, Miguel and Coniglio, Simone},
  year = {2020},
  journal = {Mechanics \& Industry},
  volume = {21},
  number = {3},
  pages = {17},
  issn = {2257-7777, 2257-7750},
  doi = {10.1051/meca/2020011},
  urldate = {2022-09-12},
  langid = {english}
}

@article{daSilva2020,
  title = {Topology Optimization of Compliant Mechanisms Considering Stress Constraints, Manufacturing Uncertainty and Geometric Nonlinearity},
  author = {{da Silva}, Gustavo Assis and Beck, Andr{\'e} Te{\'o}filo and Sigmund, Ole},
  year = {2020},
  journal = {Computer Methods in Applied Mechanics and Engineering},
  volume = {365},
  pages = {112972},
  issn = {00457825},
  doi = {10.1016/j.cma.2020.112972},
  urldate = {2020-11-10},
  langid = {english}
}

@article{daSilva2021,
  title = {Local versus Global Stress Constraint Strategies in Topology Optimization: {{A}} Comparative Study},
  shorttitle = {Local versus Global Stress Constraint Strategies in Topology Optimization},
  author = {{da Silva}, Gustavo Assis and Aage, Niels and Beck, Andr{\'e} Te{\'o}filo and Sigmund, Ole},
  year = {2021},
  journal = {International Journal for Numerical Methods in Engineering},
  volume = {122},
  number = {21},
  pages = {6003--6036},
  issn = {0029-5981, 1097-0207},
  doi = {10.1002/nme.6781},
  urldate = {2021-08-09},
  langid = {english}
}

@article{Deepak2009,
  title = {A {{Comparative Study}} of the {{Formulations}} and {{Benchmark Problems}} for the {{Topology Optimization}} of {{Compliant Mechanisms}}},
  author = {Deepak, Sangamesh R. and Dinesh, M. and Sahu, Deepak K. and Ananthasuresh, G. K.},
  year = {2009},
  journal = {Journal of Mechanisms and Robotics},
  volume = {1},
  number = {1},
  pages = {011003},
  issn = {1942-4302, 1942-4310},
  doi = {10.1115/1.2959094},
  urldate = {2020-05-19},
  langid = {english}
}

@article{DeLeon2020,
  title = {Stress-Based Topology Optimization of Compliant Mechanisms Design Using Geometrical and Material Nonlinearities},
  author = {De Leon, Daniel M. and Gon{\c c}alves, Juliano F. and {de Souza}, Carlos E.},
  year = {2020},
  journal = {Structural and Multidisciplinary Optimization},
  volume = {62},
  number = {1},
  pages = {231--248},
  issn = {1615-147X, 1615-1488},
  doi = {10.1007/s00158-019-02484-4},
  urldate = {2020-11-10},
  langid = {english}
}

@book{Gasch1989,
  title = {Strukturdynamik: {{Bd}}. 2. {{Kontinua Und Ihre Diskretisierung}}},
  author = {Gasch, Robert and Knothe, Klaus},
  year = {1989},
  volume = {2},
  publisher = {Springer},
  isbn = {3-540-50771-X}
}

@article{Golub1973b,
  title = {Some Modified Matrix Eigenvalue Problems},
  author = {Golub, Gene H.},
  year = {1973},
  journal = {Siam Review},
  volume = {15},
  number = {2},
  pages = {318--334},
  publisher = {SIAM},
  isbn = {0036-1445}
}

@article{Hasse2009,
  title = {Design of Compliant Mechanisms with Selective Compliance},
  author = {Hasse, Alexander and Campanile, Lucio Flavio},
  year = {2009},
  journal = {Smart Materials and Structures},
  volume = {18},
  number = {11},
  pages = {115016},
  issn = {0964-1726, 1361-665X},
  doi = {10.1088/0964-1726/18/11/115016},
  urldate = {2020-05-19},
  langid = {english}
}

@inproceedings{Hasse2016,
  title = {Topology {{Optimization}} of {{Compliant Mechanisms Explicitly Considering Desired Kinematics}} and {{Stiffness Constraints}}},
  booktitle = {Volume {{5A}}: 40th {{Mechanisms}} and {{Robotics Conference}}},
  author = {Hasse, Alexander},
  year = {2016},
  pages = {V05AT07A029},
  publisher = {American Society of Mechanical Engineers},
  address = {Charlotte, North Carolina, USA},
  doi = {10.1115/DETC2016-60496},
  urldate = {2020-05-19},
  isbn = {978-0-7918-5015-2},
  langid = {english}
}

@article{Joo2004,
  title = {Topological {{Synthesis}} of {{Compliant Mechanisms Using Nonlinear Beam Elements}}},
  author = {Joo, Jinyong and Kota, Sridhar},
  year = {2004},
  journal = {Mechanics Based Design of Structures and Machines},
  volume = {32},
  number = {1},
  pages = {17--38},
  issn = {1539-7734, 1539-7742},
  doi = {10.1081/SME-120026588},
  urldate = {2020-11-10},
  langid = {english}
}

@article{Kawamoto2009,
  title = {Stabilization of Geometrically Nonlinear Topology Optimization by the {{Levenberg}}--{{Marquardt}} Method},
  author = {Kawamoto, Atsushi},
  year = {2009},
  journal = {Structural and Multidisciplinary Optimization},
  volume = {37},
  number = {4},
  pages = {429--433},
  issn = {1615-147X, 1615-1488},
  doi = {10.1007/s00158-008-0236-5},
  urldate = {2022-12-07},
  langid = {english}
}

@article{Kirmse2021,
  title = {Synthesis of Compliant Mechanisms with Selective Compliance -- {{An}} Advanced Procedure},
  author = {Kirmse, Stephanie and Campanile, Lucio Flavio and Hasse, Alexander},
  year = {2021},
  journal = {Mechanism and Machine Theory},
  volume = {157},
  pages = {104184},
  issn = {0094-114X},
  doi = {10.1016/j.mechmachtheory.2020.104184}
}

@article{Klarbring2013,
  title = {Topology Optimization of Hyperelastic Bodies Including Non-Zero Prescribed Displacements},
  author = {Klarbring, Anders and Str{\"o}mberg, Niclas},
  year = {2013},
  journal = {Structural and Multidisciplinary Optimization},
  volume = {47},
  number = {1},
  pages = {37--48},
  issn = {1615-147X, 1615-1488},
  doi = {10.1007/s00158-012-0819-z},
  urldate = {2022-12-07},
  langid = {english}
}

@article{Kota1995a,
  title = {Designing Compliant Mechanisms},
  author = {Kota, Sridhar and Ananthasuresh, G. K.},
  year = {1995},
  journal = {Mechanical Engineering-CIME},
  volume = {117},
  number = {11},
  pages = {93--97},
  publisher = {American Society of Mechanical Engineers},
  isbn = {0025-6501}
}

@article{Kumar2021,
  title = {On Topology Optimization of Large Deformation Contact-Aided Shape Morphing Compliant Mechanisms},
  author = {Kumar, Prabhat and Sauer, Roger A. and Saxena, Anupam},
  year = {2021},
  journal = {Mechanism and Machine Theory},
  volume = {156},
  pages = {104135},
  issn = {0094114X},
  doi = {10.1016/j.mechmachtheory.2020.104135},
  urldate = {2020-11-10},
  langid = {english}
}

@article{Lahuerta2013,
  title = {Towards the Stabilization of the Low Density Elements in Topology Optimization with Large Deformation},
  author = {Lahuerta, Ricardo Doll and Sim{\~o}es, Eduardo T. and Campello, Eduardo M. B. and Pimenta, Paulo M. and Silva, Emilio C. N.},
  year = {2013},
  journal = {Computational Mechanics},
  volume = {52},
  number = {4},
  pages = {779--797},
  issn = {0178-7675, 1432-0924},
  doi = {10.1007/s00466-013-0843-x},
  urldate = {2022-12-07},
  langid = {english}
}

@article{Lazarov2011,
  title = {Robust Design of Large-Displacement Compliant Mechanisms},
  author = {Lazarov, B. S. and Schevenels, M. and Sigmund, O.},
  year = {2011},
  journal = {Mechanical Sciences},
  volume = {2},
  number = {2},
  pages = {175--182},
  issn = {2191-916X},
  doi = {10.5194/ms-2-175-2011},
  urldate = {2020-11-10},
  langid = {english}
}

@article{Li2014a,
  title = {Reliability-{{Based Topology Optimization}} of {{Compliant Mechanisms}} with {{Geometrically Nonlinearity}}},
  author = {Li, Zhao Kun and Bian, Hua Mei and Shi, Li Juan and Niu, Xiao Tie},
  year = {2014},
  journal = {Applied Mechanics and Materials},
  volume = {556},
  pages = {4422--4434},
  issn = {1662-7482},
  doi = {10.4028/www.scientific.net/AMM.556-562.4422},
  urldate = {2020-11-10},
  langid = {english}
}

@article{Liu2017a,
  title = {Design of {{Large-Displacement Compliant Mechanisms}} by {{Topology Optimization Incorporating Modified Additive Hyperelasticity Technique}}},
  author = {Liu, Liying and Xing, Jian and Yang, Qingwei and Luo, Yangjun},
  year = {2017},
  journal = {Mathematical Problems in Engineering},
  volume = {2017},
  pages = {1--11},
  issn = {1024-123X, 1563-5147},
  doi = {10.1155/2017/4679746},
  urldate = {2020-04-29},
  langid = {english}
}

@article{Luo2015,
  title = {Topology Optimization of Geometrically Nonlinear Structures Based on an Additive Hyperelasticity Technique},
  author = {Luo, Yangjun and Wang, Michael Yu and Kang, Zhan},
  year = {2015},
  journal = {Computer Methods in Applied Mechanics and Engineering},
  volume = {286},
  pages = {422--441},
  issn = {00457825},
  doi = {10.1016/j.cma.2014.12.023},
  urldate = {2022-12-07},
  langid = {english}
}

@article{Mankame2007,
  title = {Synthesis of Contact-Aided Compliant Mechanisms for Non-Smooth Path Generation},
  author = {Mankame, N. D. and Ananthasuresh, G. K.},
  year = {2007},
  journal = {International Journal for Numerical Methods in Engineering},
  volume = {69},
  number = {12},
  pages = {2564--2605},
  issn = {00295981, 10970207},
  doi = {10.1002/nme.1861},
  urldate = {2022-09-12},
  langid = {english}
}

@article{Maute2003,
  title = {Reliability-Based Design of {{MEMS}} Mechanisms by Topology Optimization},
  author = {Maute, Kurt and Frangopol, Dan M.},
  year = {2003},
  journal = {Computers \& Structures},
  volume = {81},
  number = {8-11},
  pages = {813--824},
  issn = {00457949},
  doi = {10.1016/S0045-7949(03)00008-7},
  urldate = {2020-04-29},
  langid = {english}
}

@article{Nelder1965,
  title = {A {{Simplex Method}} for {{Function Minimization}}},
  author = {Nelder, J. A. and Mead, R.},
  year = {1965},
  journal = {The Computer Journal},
  volume = {7},
  number = {4},
  pages = {308--313},
  issn = {0010-4620},
  doi = {10.1093/comjnl/7.4.308},
  urldate = {2021-01-26}
}

@article{Pedersen2001,
  title = {Topology Synthesis of Large-Displacement Compliant Mechanisms},
  author = {Pedersen, Claus B. W. and Buhl, Thomas and Sigmund, Ole},
  year = {2001},
  journal = {International Journal for Numerical Methods in Engineering},
  volume = {50},
  number = {12},
  pages = {2683--2705},
  issn = {0029-5981, 1097-0207},
  doi = {10.1002/nme.148},
  urldate = {2020-04-29},
  langid = {english}
}

@article{Saxena2001a,
  title = {Topology {{Synthesis}} of {{Compliant Mechanisms}} for {{Nonlinear Force-Deflection}} and {{Curved Path Specifications}}},
  author = {Saxena, A. and Ananthasuresh, G. K.},
  year = {2001},
  month = mar,
  journal = {Journal of Mechanical Design},
  volume = {123},
  number = {1},
  pages = {33--42},
  issn = {1050-0472, 1528-9001},
  doi = {10.1115/1.1333096},
  urldate = {2020-06-17},
  langid = {english}
}

@article{Seltmann2022,
  title = {Topology-Optimization Based Design of Multi-Degree-of-Freedom Compliant Mechanisms (Mechanisms with Multiple Pseudo-Mobility)},
  author = {Seltmann, Stephanie and Campanile, Lucio Flavio and Hasse, Alexander},
  year = {2022},
  journal = {Journal of Intelligent Material Systems and Structures},
  volume = {34},
  number = {5},
  pages = {609--628},
  doi = {10.1177/1045389X221111535},
  urldate = {2022-02-03}
}

@article{Seltmann2023,
  title = {Topology Optimization of Compliant Mechanisms with Distributed Compliance (Hinge-Free Compliant Mechanisms) by Using Stiffness and Adaptive Volume Constraints Instead of Stress Constraints},
  author = {Seltmann, Stephanie and Hasse, Alexander},
  year = {2023},
  journal = {Mechanism and Machine Theory},
  volume = {180},
  number = {11},
  pages = {105133},
  issn = {0094114X},
  doi = {10.1016/j.mechmachtheory.2022.105133},
  urldate = {2022-12-07},
  langid = {english}
}

@incollection{Simo1998a,
  title = {A Simple Model of Hyperelastic Response},
  booktitle = {Computational {{Inelasticity}}},
  author = {Simo, Juan C. and Hughes, Thomas JR},
  year = {1998},
  edition = {1},
  volume = {7},
  pages = {258--259},
  publisher = {Springer Science \& Business Media},
  isbn = {978-0-387-97520-7}
}

@article{Wang2014c,
  title = {Interpolation Scheme for Fictitious Domain Techniques and Topology Optimization of Finite Strain Elastic Problems},
  author = {Wang, Fengwen and Lazarov, Boyan Stefanov and Sigmund, Ole and Jensen, Jakob S{\o}ndergaard},
  year = {2014},
  month = jul,
  journal = {Computer Methods in Applied Mechanics and Engineering},
  volume = {276},
  pages = {453--472},
  issn = {00457825},
  doi = {10.1016/j.cma.2014.03.021},
  urldate = {2022-12-07},
  langid = {english}
}

@inproceedings{XavierLeitao2019,
  title = {Topology {{Optimization}} of {{Geometrically Nonlinear Structures Considering}} an {{Energy Interpolation Scheme}}},
  booktitle = {Proceedings of the 25th {{International Congress}} of {{Mechanical Engineering}}},
  author = {Xavier Leit{\~a}o, Andr{\'e} and Pereira, Anderson},
  year = {2019},
  publisher = {ABCM},
  doi = {10.26678/ABCM.COBEM2019.COB2019-1580},
  urldate = {2022-12-07},
  langid = {english}
}

@article{Zhu2020,
  title = {Design of Compliant Mechanisms Using Continuum Topology Optimization: {{A}} Review},
  shorttitle = {Design of Compliant Mechanisms Using Continuum Topology Optimization},
  author = {Zhu, Benliang and Zhang, Xianmin and Zhang, Hongchuan and Liang, Junwen and Zang, Haoyan and Li, Hai and Wang, Rixin},
  year = {2020},
  journal = {Mechanism and Machine Theory},
  volume = {143},
  pages = {103622},
  issn = {0094114X},
  doi = {10.1016/j.mechmachtheory.2019.103622},
  urldate = {2020-09-30},
  langid = {english}
}

\end{document}